\documentclass{amsart}

\usepackage{physics}
\usepackage{subcaption}
\usepackage{graphicx}  % for "\includegraphics" macro
\usepackage{longtable}

\theoremstyle{definition}

\theoremstyle{remark}

\numberwithin{equation}{section}

\begin{document}

\title{Green's Function of the Screened Poisson's Equation on the Sphere}

%    Information for first author
\author{Ramy Tanios}
\address{American University of Beirut, Beirut, Lebanon}
\curraddr{xxx}
\email{rgt09@mail.aub.edu}
%\thanks{The first author was supported in part by NSF Grant \#000000.}

%    Information for second author
\author{Samah El Mohtar}
\address{American University of Beirut, Beirut, Lebanon}
\curraddr{King Abdullah University of Science and Technology, Thuwal 23955, Saudi Arabia}
\email{Samah.Mohtar@kaust.edu.sa}

\author{Omar Knio}
\address{Duke University, Durham, NC 27708, USA}
\curraddr{King Abdullah University of Science and Technology, Thuwal 23955, Saudi Arabia}
\email{Omar.Knio@kaust.edu.sa}

\author{Issam Lakkis}
\address{American University of Beirut, Beirut, Lebanon}
\email{issam.lakkis@aub.du.lb}

%    General info
%\subjclass[2000]{Primary 54C40, 14E20; Secondary 46E25, 20C20}

\date{\today}

%\dedicatory{This paper is dedicated to our advisors.}


\begin{abstract}
In geophysical fluid dynamics, the screened Poisson equation appears in the shallow-water, quasi geostrophic equations. Recently, many attempts have been made to solve those equations on the sphere using different numerical methods.  These include vortex methods, which solve a Poisson equation to compute the streamfunction from the (relative) vorticity. Alternatively, the streamfunction can be computed directly from potential vorticity (PV), which would offer the possibility of constructing more attractive vortex methods because PV is conserved along material trajectories in the inviscid case. On the spherical shell, however, the screened Poisson equation does not admit a known Green's function, which limits the extension of such approaches to the case of a sphere. In this paper, we derive an expression of Green's function for the screened Poisson equation on the spherical shell in series form and in integral form. A proof of convergence of the series representation is then given. As the series is slowly convergent, a robust and efficient approximation is obtained using a split form which isolates the singular behavior.  The solutions are illustrated and analyzed for different values of the screening constant.
\end{abstract}

\maketitle

\section{Introduction}
\label{sec:intro}
In geophysical fluid dynamics, the screened Poisson equation arises in the shallow-water, quasi-geostrophic, potential vorticity equation~\cite{Pedlosky1987,Vallis2006}, for a finite Rossby radius of deformation (baroclinic case). The equation relates the 
streamfunction of the geostrophic flow to the potential vorticity, where the ``screening'' is the inverse of the Rossby radius of deformation. 
For an infinite Rossby radius of deformation (barotropic case), the screened Poisson equation reduces to the Poisson equation, whose Green's function is known on the spherical shell~\cite{Bogomolov1977,KimuraOkamoto1987}. Recently, many attempts using Lagrangian methods have been made to (numerically) solve the shallow-water quasi-geostrophic potential vorticity equation on the spherical shell. 
For instance, Bosler et al.~\cite{Bosler2014} solved the barotropic vorticity (BVE) equation (infinite Rossby radius of deformation) using a Lagrangian particle/panel method. The flow field was computed from the position of particles carrying (relative) vorticity, and advecting with
a velocity expressed in terms of the Biot-Savat law. However, the method did not take advantage of the conservation of potential vorticity, i.e., (relative) vorticity carried by each particle had to be updated at the new particle positions, thus requiring an additional computational cost. 
In \cite{MohammadianMarshall2010}, Mohammadian \& Marshall used a vortex-in-cell (VIC) method, in which particles carried (relative) vorticity. The flow field was obtained from the streamfunction, which was computed from the vorticity by inverting a Poisson equation on an underlying Eulerian grid.  

Allowing particles to carry potential vorticity enables taking advantage of the fact that potential vorticity is materially conserved in the 
inviscid limit. Consequently, in this case advecting/transporting particles along flow trajectories would avoid the need to integrate an
evolution equation for their strengths, provided that the flow field can be immediately computed from the particle distribution.  To this end,
we focus on this work on deriving expressions of Green's function for the screened Poisson equation on the spherical shell.

Specifically, in section~\ref{sec:gf-der}, we apply a spectral decomposition of the Laplace-Beltrami operator to arrive at a series 
representation of the Green's function.  The convergence properties of this representation are then analyzed in section~\ref{sec:gf-conv},
and a computational strategy for evaluating the series is outlined in section~\ref{sec:numer}.  In section~\ref{sec:int-form}, an alternative, 
integral form of the Green's function is constructed.  Implementation of the series and integral solutions is then illustrated in 
section~\ref{sec:res}, in light of results obtained for representative test cases.  Concluding remarks are given in section~\ref{sec:conc}.

%\begin{itemize}
%\item Derive the Green's function of the Screened Poisson's Equation on the Sphere.
%\item Prove the convergence of the Green's funcion.
%\item Derive the integral form.
%\item Numerically Evaluate the series form of the Green's function.
%\end{itemize}
%

%%%%%%%%%%%%%%%%%%%%%%%%%%%%%%%%%%%%%%%%%%%%%%%%%%%%
%\newpage
%%%%%%%%%%%%%%%%%%%%%%%%%%%%%%%%%%%%%%%%%%%%%%%%%%%%%
\section{Derivation of Green's function}
\label{sec:gf-der}

Let $\Omega = \{(\rho,\theta,\varphi) \in \mathbb{R}^+\times[0,2\pi]\times[0,\pi]  \  /  \ \rho = R\}$. 
Consider the screened Poisson's equation on $\Omega$: 
\begin{equation}
\nabla^2_s \psi(\theta,\varphi) - \frac{1}{L_d^2}\psi(\theta,\varphi)= f(\theta,\varphi)
\label{eq:SP}
\end{equation}
where $L_d \in \mathbb{R}^{+}$ is the Rossby radius of deformation and $\nabla^2_s$ is the \textbf{Laplace-Beltrami} operator on the sphere of radius $R$,
\begin{equation}
\nabla^2_s = \frac{1}{R^2}\frac{1}{\sin^2 \theta }\pdv[2]{}{\varphi} + \frac{1}{R^2}\frac{1}{\sin\theta}\pdv[]{}{\theta} (\sin\theta \pdv[]{}{\theta}).
\end{equation}

The spherical harmonics $Y_{l,m}(\theta,\varphi)$ form a complete basis set of the Hilbert space of all square-integrable functions,
${\mathcal H} = \left\{ f: \Omega \to \mathbb{R} / \int_{\Omega} f^2 < \infty \right\}$. Thus every function in ${\mathcal H}$ can be decomposed
in terms of the mean-square convergent sum:
\begin{equation}
f(\theta,\varphi) = \sum_{l=0}^{\infty} \sum_{m=-l}^l f_{lm}(r)Y_{l,m}(\theta,\varphi),
\label{eq:f-harmonic}
\end{equation}
and the solution, $\psi$, of (\ref{eq:SP}) can be expressed as:
\begin{equation}
\psi(\theta,\varphi) = \sum_{l=0}^{\infty} \sum_{m=-l}^l u_{lm}(r)Y_{l,m}(\theta,\varphi).
\end{equation}
Using the $L_2$ inner product:
\begin{equation}
 <h(\theta,\varphi),k(\theta,\varphi)>  \equiv \int\limits_{\Omega} h(\theta,\varphi)k(\theta,\varphi) dS  \quad (h,k) \in {\mathcal H}^2
\end{equation}
the coefficients in (\ref{eq:f-harmonic}) are given by:
\begin{equation}
f_{lm} = \frac{1}{R^2}\int\limits_{\Omega} Y_{l,m}(\theta,\varphi) f(\theta,\varphi) dS  = \int_{\theta=0}^{\pi} \int_{\varphi=0}^{2\pi}Y_{l,m}(\theta,\varphi)f(\theta,\varphi) \sin\theta \,d\theta\,d\varphi.
\end{equation}
Because the spherical harmonics are eigenfunctions of $\nabla^2_s|_{R=1}$~\cite{SphericalFunctions}, that is 
\begin{equation}
\left. \nabla^2_s \right| _{R=1} Y_{l,m}  = -l(l+1) Y_{l,m}
\end{equation}
we have 
\begin{equation}
\nabla^2_s Y_{l,m}  =\frac{1}{R^2} \left. \nabla^2_{s} \right| _{R=1}  Y_{l,m}=\frac{-l(l+1)}{R^2}Y_{l,m}
\end{equation}
Now we write (\ref{eq:SP}) as:
\begin{equation}
\sum_{l=0}^{\infty} \sum_{m=-l}^l \Big[ u_{lm} \frac{(-l)(l+1)}{R^2} Y_{l,m}(\theta,\varphi)-\frac{1}{L_d^2}u_{lm}(r)Y_{l,m}(\theta,\varphi)\Big] = \sum_{l=0}^{\infty} \sum_{m=-l}^l f_{lm}Y_{l,m}(\theta,\varphi) .
\end{equation}
From the orthogonality of the basis, we obtain
\begin{equation}
u_{lm} = \frac{ f_{lm}}{ \frac{(-l)(l+1)}{R^2}-\frac{1}{L_d^2}}= \frac{ \int_{\theta'=0}^{\pi} \int_{\varphi'=0}^{2\pi}Y_{l,m}(\theta',\varphi')f(\theta',\varphi') \sin\theta' \,d\theta'\,d\varphi'
}{ \frac{(-l)(l+1)}{R^2}-\frac{1}{L_d^2}}
\end{equation}
and 
\begin{align}
\psi(\theta,\varphi) &= \sum_{l=0}^{\infty} \sum_{m=-l}^l \frac{  f_{lm}}{ \frac{(-l)(l+1)}{R^2}-\frac{1}{L_d^2}}Y_{l,m}(\theta,\varphi) \\&=\sum_{l=0}^{\infty} \sum_{m=-l}^l \frac{ \int_{\theta'=0}^{\pi} \int_{\varphi'=0}^{2\pi}Y_{l,m}(\theta',\varphi')f(\theta',\varphi') \sin\theta' \,d\theta'\,d\varphi'
}{ \frac{(-l)(l+1)}{R^2}-\frac{1}{L_d^2}}Y_{l,m}(\theta,\varphi) \\ &= 
\int_{\theta'=0}^{\pi} \int_{\varphi'=0}^{2\pi} \sum_{l=0}^{\infty} \sum_{m=-l}^l \frac{ Y_{l,m}(\theta,\varphi)Y_{l,m}(\theta',\varphi')}{(-l)(l+1)- \frac{R^2}{L_d^2}}f(\theta',\varphi') R^2 \sin\theta' \,d\theta'\,d\varphi'
\end{align}
%Knowing that the bounds of the summations and the bounds of the integrals are independent. \\
%\vspace{3mm}

Because $\psi = G * f$, we have:
\begin{align}
G((R,\theta,\varphi),(R,\theta',\varphi')) & =\sum_{l=0}^{\infty} \sum_{m=-l}^l \frac{ Y_{l,m}(\theta',\varphi')Y_{l,m}(\theta,\varphi)}{(-l)(l+1)-\frac{R^2}{L_d^2}} \\
 & = \sum_{l=0}^{\infty}\frac{1}{(-l)(l+1)-\frac{R^2}{L_d^2}}\sum_{m=-l}^l  Y_{l,m}(\theta',\varphi')Y_{l,m}(\theta,\varphi)
\end{align}
Using (i) the spherical harmonics addition theorem~\cite{SphericalFunctions}:
\begin{equation}
\forall (R,\theta,\varphi),(R,\theta',\varphi') \in \Omega,  \ \frac{4\pi}{2l+1} \sum_{m=-l}^{l}Y_{l,m}(\theta,\varphi)Y_{l,m}(\theta',\varphi') = P_l(\cos\gamma)
\end{equation}
where $\gamma$ is the angle at the center between  $(R,\theta,\varphi)$ and $(R,\theta',\varphi')$, and (ii) the identity $\cos\gamma =\cos\theta\cos\theta'+\sin\theta\sin\theta'\cos(\varphi-\varphi')$, we obtain
\begin{align}
G((R,\theta,\varphi),(R,\theta',\varphi')) &= \sum_{l=0}^{\infty}\frac{1}{(-l)(l+1)-\frac{R^2}{L_d^2}} \frac{2l+1}{4\pi}P_l(\cos(\gamma))\\ &= \boxed{\frac{-1}{4\pi}\sum_{l=0}^{\infty}\frac{2l+1}{l(l+1)+\frac{R^2}{L_d^2}}P_l(\cos(\gamma))}
\label{eq:G-series}
\end{align}

%%%%%%%%%%%%%%%%%%%%%%%%%%%%%%%%%%%%%%%%%%%%%%%%%%%%
%\newpage
\section{Convergence of the series representation}
\label{sec:gf-conv}

In this section, we briefly examine properties of the Green's function series representation.
To this end, we study the behavior of:
\begin{equation}
f\left(w,\gamma\right)=\sum_{l\geq0}\frac{\left(2l+1\right)P_{l}\left(\cos\left(\gamma\right)\right)}{l\left(l+1\right)+w},\,w\in\mathbb{R}^{+}.
\end{equation}

For $\cos(\gamma)=1$, we have $P_{l} \left(\cos\left(\gamma\right)\right)=1,\ \forall l \in\mathbb{N}$, consequently
\begin{equation}
\sum_{l\geq0}\frac{2l+1}{l\left(l+1\right)+w}\sim\sum_{l\geq0}\frac{2}{l+w} , 
\end{equation}
and so the series diverges like the harmonic series. 

For $\cos(\gamma) \neq 1$, we have
\begin{equation}
f\left(w,\gamma\right)=2\sum_{l\geq0}\frac{lP_{l}\left(\cos\left(\gamma\right)\right)}{l\left(l+1\right)+w}+\sum_{l\geq0}\frac{P_{l}\left(\cos\left(\gamma\right)\right)}{l\left(l+1\right)+w}.
\label{eq:f2series}
\end{equation}
The second series is absolutely convergent because
\begin{equation}
\sum_{l\geq0}\left|\frac{P_{l}\left(\cos\left(\gamma\right)\right)}{l\left(l+1\right)+w}\right|\leq\sum_{l\geq0}\frac{1}{l^{2}+w}<\infty,\,\forall w\in\mathbb{R}^{+}
\end{equation}
Now consider the first series,
let $A_l = \frac{l}{l(l+1)+w}$ and $B_l = P_l(\cos(\gamma))$.  Clearly, $A_l \geq 0$ and $\lim_{l\to\infty} A_l = 0$.
Hence, there exists  $ l^* \in \mathbb{N}$ such that  $A_{l+1}\leq A_l$ for all  $l\ \geq l^*$.
We now rewrite the first series on the right hand side of (\ref{eq:f2series}) as: 
\begin{equation}
\underbrace{\sum_{l=0}^{l^*}\frac{lP_{l}\left(\cos\left(\gamma\right)\right)}{l\left(l+1\right)+w}}_{\text{finite sum}} + \sum_{k=0}^{\infty}\frac{(l^*+k+1)P_{(l^*+k+1)}\left(\cos\left(\gamma\right)\right)}{(l^*+k+1)\left((l^*+k+2)\right)+w}.
\end{equation}
Let us show that the second member is convergent using the Dirichlet Test. We have already shown that $A_l \to 0$ and that it is monotonically decreasing for $l \geq l^*$. 
It remains to show that  there exists $M \geq 0$ such that  $ \left| \sum_{l=0}^N P_l(\cos(\gamma)) \right| \leq M$ for all $N$.

Note that if $\cos(\gamma)=-1.0$, the $P_l\left( \cos(\gamma) \right) = (-1)^l P_l(1) = (-1)^l$, therefore $\sum_{l=0}^{N} P_l \left( \cos(\gamma) \right) = \pm 1$, which is bounded.

We may thus assume that $ \left| \cos(\gamma) \right| < 1$.  We make use of (i) the Legendre generating function:
\begin{equation}
\sum_{l=0}^{\infty} u^l P_l(y) = \frac{1}{\sqrt{u^2 - 2yu+1}},
\end{equation}
with $y=\cos(\gamma)$, and (ii) the binomial series expansion $(1+x)^{1/2}$ with $x = u^2 - 2y$.

For $0<\cos(\gamma)<1$, we let $u=1$, so we have  $x=1-2\cos(\gamma) <1$.  Consequently, the binomial series converges and the 
Legendre sum is bounded.  For $-1<\cos(\gamma)<0$, we use the fact that $P_l(-\cos(\gamma)) = (-1)^l P_l(\cos(\gamma))$, let $u = -1$ and $x = 1-2|\cos(\gamma)|.$   Following the same argument as before, we conclude that the sum is bounded.  Finally, if $\cos(\gamma)=0$, we may set 
$u \pm 1$, which leads $x = 1$, and the same conclusion as before.

Consequently, by the Dirichlet test, the series $\sum_{l=0}^{\infty}\frac{lP_{l}\left(\cos\left(\gamma\right)\right)}{l\left(l+1\right)+w}$ converges when 
$\cos(\gamma) < 1$.

%%%%%%%%%%%%%%%%%%%%%%%%%%%%%%%%%%%%%%%%%%%%%
%\newpage
\section{Numerical approximation of the series representation}
\label{sec:numer}

In the tests below, we assess two approaches for estimating the Green's function based on its
series representation.  The first is a straightforward approach based on truncating (\ref{eq:G-series}) at
a suitably large index, $l'$, namely through:
\begin{equation}
G((R,\theta,\varphi),(R,\theta',\varphi')) \approx \frac{-1}{4\pi}\sum_{l=0}^{l'}\frac{2l+1}{(l)(l+1)+\frac{R^2}{L_d^2}}P_l(\cos(\gamma)) .
\label{eq:G-trunc}
\end{equation}
We refer to (\ref{eq:G-trunc}) as the truncated approximation.

A second, alternative approach is developed based on first splitting (\ref{eq:G-series}) according to:
\begin{equation}
G((R,\theta,\varphi),(R,\theta',\varphi')) = \frac{-1}{4\pi R^2}\sum_{l=0}^{l'-1}\frac{2l+1}{\frac{l(l+1)}{R^2}+\frac{1}{L_d^2}}P_l(\cos(\gamma)) +
\frac{-1}{4\pi R^2}\sum_{l=l'}^{\infty}\frac{2l+1}{\frac{l(l+1)}{R^2}+\frac{1}{L_d^2}}P_l(\cos(\gamma)) .
\label{eq:G-alt1}
\end{equation}
Selecting $l'$ such that $\frac{l'(l'+1)}{R^2} >> \frac{1}{L_d^2}$, we may approximate $G$ according to:
\begin{equation}
G((R,\theta,\varphi),(R,\theta',\varphi')) \approx \frac{-1}{4\pi R^2}\sum_{l=0}^{l'-1}\frac{2l+1}{\frac{l(l+1)}{R^2}+\frac{1}{L_d^2}}P_l(\cos(\gamma)) +
\frac{-1}{4\pi R^2}\sum_{l=l'}^{\infty}\frac{2l+1}{\frac{l(l+1)}{R^2}}P_l(\cos(\gamma)) .
\label{eq:G-alt2}
\end{equation}
Let
\begin{equation} 
G^*((\theta,\varphi),(\theta',\varphi')) = \frac{1}{4\pi}\log(\frac{\text{e}}{2}(1-\cos(\gamma)))
\end{equation}
denote the Green's function of the Poisson equation on the sphere,
\begin{equation}
\nabla_s^2\psi(\theta,\varphi) = f(\theta,\varphi).
\end{equation}
$G^*$ can be expressed in series form as:
\begin{equation}
G^*((\theta,\varphi),(\theta',\varphi')) = \frac{-1}{4\pi}\sum_{l=1}^{\infty}\frac{2l+1}{(l)(l+1)}P_l(\cos(\gamma))
\label{eq:Gstar}
\end{equation}
Inserting (\ref{eq:Gstar}) into (\ref{eq:G-alt2}) and rearranging we finally obtain:
\begin{align}
G((R,\theta,\varphi),(R,\theta',\varphi')) &\approx -\frac{L_d^2}{4\pi R^2}-\frac{1}{4\pi R^2}\sum_{l=1}^{l'-1} \left[\frac{2l+1}{\frac{(l)(l+1)}{R^2} + \frac{1}{L_d^2}} - \frac{2l+1}{\frac{(l)(l+1)}{R^2}}
\right] P_l(\cos(\gamma))  \nonumber \\&+ \frac{1}{4\pi} \log(\frac{e}{2}(1-\cos(\gamma))) ,
\label{eq:G-split}
\end{align}
which we refer to as the split approximation.
Because the Rossby radius, $L_d$, defines a distance on the circumference of the sphere, and the distance between the source and target on the sphere is $R \gamma$, we introduce the characteristic 
angle of the problem $\gamma^* \equiv L_d/R$. The Green's function and its split sum approximation are then expressed in terms of $\gamma^*$ as:
\begin{align}
 G(\gamma,\gamma^*) \approx G_{l'}(\gamma,\gamma^*) & = -\frac{\gamma^{*2}}{4\pi }-\frac{1}{4\pi}\sum_{l=1}^{l'-1} \left[\frac{2l+1}{(l)(l+1)+ \frac{1}{\gamma^{*2}} } - \frac{2l+1}{(l)(l+1)} \right] P_l(\cos(\gamma))  \nonumber \\&+ \frac{1}{4\pi} \log(\frac{e}{2}(1-\cos(\gamma))) .
\label{eq:G-split1}
\end{align}
As further discussed below, with the same truncation index, $l'$, the split approximation leads to estimates exhibiting 
appreciably smaller relative errors than straightforward truncation.  It also exhibits faster and more robust convergence
as $l'$ increases.

%%%%%%%%%%%%%%%%%%%%%%%%%%%%%%%%%%%%%%%%%%%%%%%%%%%%
%\newpage
\section{Integral form}
\label{sec:int-form}
In this section, we exploit the series solution (\ref{eq:G-series}) to derive an alternative, integral form of the Green's function.
To this end, we make use of the following identity,
\begin{equation}
\frac{1}{l+R} = \int_{0}^{+\infty} e^{-z(l+R)} \,dz, \quad {\rm Re}\left( l+R \right) > 0 , 
\label{eq:identity}
\end{equation}
and factor term $l(l+1) + \frac{R^2}{L_d^2}$ as $(l-S_1)(l-S_2)$ where $S_1,S_2 \in \mathbb{C}$.  

We then perform the partial fraction expansion,
$$\frac{2l+1}{l(l+1) +  \frac{R^2}{L_d^2}} = \frac{s_1}{l-S_1} + \frac{s_2}{l-S_2},$$
where $s_1,s_2 \in \mathbb{C}$.  Note that $S_1$ and $S_2$ are complex conjugates that are independent of $l$, and 
so are $s_1$ and $s_2$.

Next, we apply (\ref{eq:identity}) to re-express the fractions $s_1/(l-S_1)$ and $s_2/(l-S_2)$ respectively according to:
$$\frac{s_1}{l-S_1} = s_1 \int_{0}^{+\infty} e^{-z(l-S_1)} \,dz ,$$ 
and 
$$\frac{s_2}{l-S_2} = s_2 \int_{0}^{+\infty} e^{-z(l-S_2)} \,dz.$$
Substituting these representations into (\ref{eq:G-series}), we get:
\begin{align}
G(R,\theta,\phi) &= \frac{-1}{4\pi}\sum_{l=0}^{\infty} (\int_{0}^{+\infty} s_1e^{-z(l-S_1)} + s_2 e^{-z(l-S_2)} \,dz)P_l(\cos(\gamma)) \\
&= \frac{-1}{4\pi}\sum_{l=0}^{\infty} (\int_{0}^{+\infty} e^{-zl}[s_1e^{zS_1} + s_2 e^{zS_2}] \,dz)P_l(\cos(\gamma)) \\
&= \frac{-1}{4\pi}\int_{0}^{+\infty} (s_1e^{zS_1} + s_2 e^{zS_2}) \sum_{l=0}^{\infty} (e^{-z})^l P_l(\cos(\gamma)) \,dz . \label{eq:leg-sum}
\end{align}
Finally, we use the Legendre generating function to re-express the summation in (\ref{eq:leg-sum}), which results in:
\begin{equation}
G(R,\theta,\phi)  = \frac{-1}{4\pi} \int_{0}^{+\infty} \frac{(s_1e^{zS_1} + s_2 e^{zS_2})}{\sqrt{e^{-2z} - 2e^{-z} \cos(\gamma) + 1}} dz .
\label{eq:G-int1}
\end{equation}

A more convenient form of (\ref{eq:G-int1}) can be obtained by substituting the values of $S_1$, $S_2$, $s_1$ and $s_2$, namely
$s_1 = s_2 = 1$,  and
$$
S_1 = -\frac{1}{2} + i \beta , \quad S_1 = -\frac{1}{2} - i \beta ,
$$
where
$$
\beta \equiv \sqrt{  \frac{R^2}{L_d^2} - \frac{1}{4} } .
$$
Performing the substitution and rearranging, we obtain:
\begin{equation}
G(R,\theta,\phi)  = \frac{-1}{2\pi} \int_{0}^{+\infty} \frac{ e^{-z/2} \cos(\beta z)}{\sqrt{e^{-2z} - 2e^{-z} \cos(\gamma) + 1}} dz .
\label{eq:G-int-form}
\end{equation}
Below, we use a quadrature approximation of (\ref{eq:G-int-form}) to verify results obtained using the series representation.

Note that for particular values of $\cos(\gamma)$, the integral in (\ref{eq:G-int-form}) can be evaluated analytically with the
result expressed in terms of elementary functions.  Specifically, for $\cos(\gamma) = -1$, we have (see~\cite{GR}, article 3.981):
\begin{equation}
\left. G \right|_{\cos(\gamma) = -1} = \frac{-1}{4 \cosh(\pi\beta)} 
\label{eqCOSM1}
\end{equation}
whereas for $\cos(\gamma) = 0$ we obtain (see~\cite{GR}, article 3.985):
\begin{equation}
\left. G \right|_{\cos(\gamma) = 0} = \frac{-1}{8 \pi^{3/2}} \Gamma \left(  \frac{1}{4} + i \frac{\beta}{2}    \right)  \Gamma \left(  \frac{1}{4} - i \frac{\beta}{2}    \right) 
\label{eqCOS0}
\end{equation}
where $\Gamma$ denotes the gamma function.  In particular, we use the results to verify both our quadrature and series approximations.

Also note the integral representation in (\ref{eq:G-int-form}) can alternatively be expressed in terms of associated Legendre function,
namely according to:
\begin{equation}
G(R,\theta,\phi)  = - \frac{1}{4 \cosh(\pi \beta)} P_{-\frac{1}{2} + i \beta}\left( \cos(\pi - \gamma) \right)
\label{eq:alf}
\end{equation}
Because the zeros of $P_{-\frac{1}{2} + i \beta}( z )$ are all real and greater than unity (see~\cite{GR}, article 8.784), we conclude
that $G(R,\theta,\phi)$ is negative in the entire range $0 < \gamma \leq \pi$ whereas it diverges to $-\infty$ as $\gamma \to 0$.

%%%%%%%%%%%%%%%%%%%%%%%%%%%%%%%%%%%%%%%%%%%%%%%

%\newpage
\section{Results}
\label{sec:res}

In this section, we first compare the split sum to the direct sum in terms of their behavior as a function of the number of terms, for the case $L_d = 1000 \ \text{km}$. The computations are carried out in Fortran using quadruple precision for real numbers and double precision for integers.
Second, we show the absolute error versus the number of terms for the split sum approximation, for  $L_d=100 \ \text{km}$ and $L_d=1000 \ \text{km}.$ These two values 
are selected to represent the low and high values of the Rossby radius of deformation in the ocean and the atmosphere\cite{Gill82,chelton98}. 
Third, we compare, for $L_d = 1000 \ \text{km}$, the values of the Green's function computed using the split sum approximation to those computed using the high precision Numerical Integration Polyalgorithm of Maple$^\text{TM}$.
  Fourth, we present plots, for different values of $L_d$, of the Green's function computed using the split sum versus the angle.  Finally, tables of the Green's function versus the angle for different values of $L_d$ are also presented in the appendix.

\begin{figure}[h]
\centering
\begin{subfigure}{0.5\textwidth}
  \centering
  \includegraphics[width=1\linewidth]{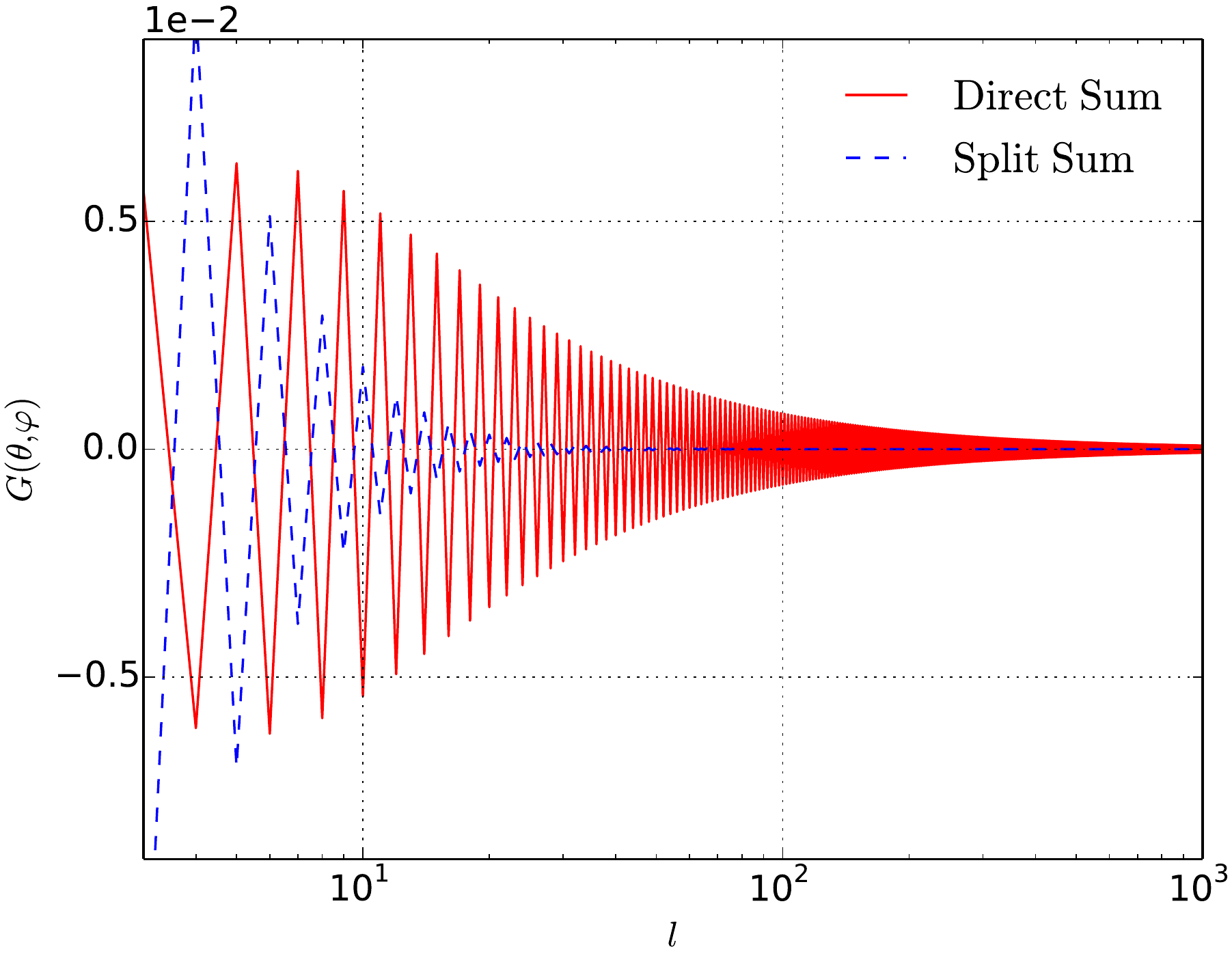}
  \caption{$\cos(\gamma)=-1$}
  \label{fig:sub1}
\end{subfigure}%
\begin{subfigure}{0.5\textwidth}
  \centering
  \includegraphics[width=1\linewidth]{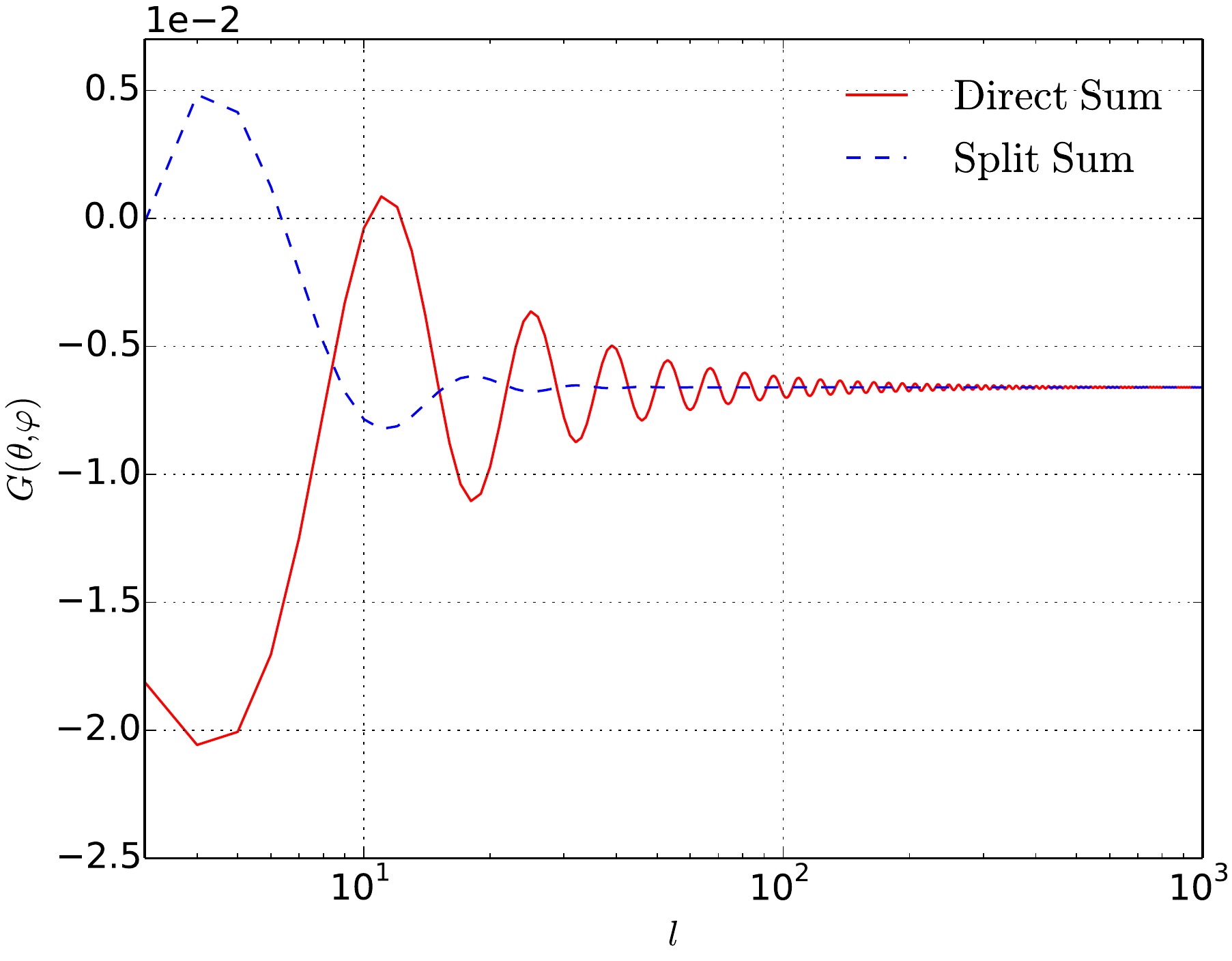}
  \caption{$\cos(\gamma)=0.9$}
  \label{fig:sub2}
\end{subfigure}
\caption{$G(\theta,\varphi)$ versus the number of terms for (a) $\cos(\gamma)=-1$, (b) $\cos(\gamma)=0.9$. Curves are generated for $L_d = 1000 \ \text{km}$ using the direct and split sum as indicated.}
\label{fig:test}
\end{figure}

Figure~\ref{fig:test} depicts the estimates obtained using the split sum and the direct sum as a function of the number of terms retained in the corresponding expansions. 
Two cases are considered, namely $\cos(\gamma)=-1.0$ and $\cos(\gamma) = 0.9$. The first case corresponds to a maximum separation between the target and the source. In the second case, the angle separating the target from the source is small. In both cases $L_d = 1000 \ \text{km}$. We observe from Figures \ref{fig:sub1} and \ref{fig:sub2} that the split sum converges faster than the direct sum. One can also observe that the rate of convergence is slower for $\cos(\gamma)=-1.0$. In fact, for a given value of $L_d$, the rate of convergence of the split sum is the slowest when $\cos(\gamma)=-1.0$. This is  because the Legendre polynomial in the split sum approximation of equation (\ref{eq:G-split}) switches sign every term, i.e. $P_l(-1)P_{l+1}(-1)<0$.

\begin{figure}[h]
  \centering
    \includegraphics[width=0.8\textwidth]{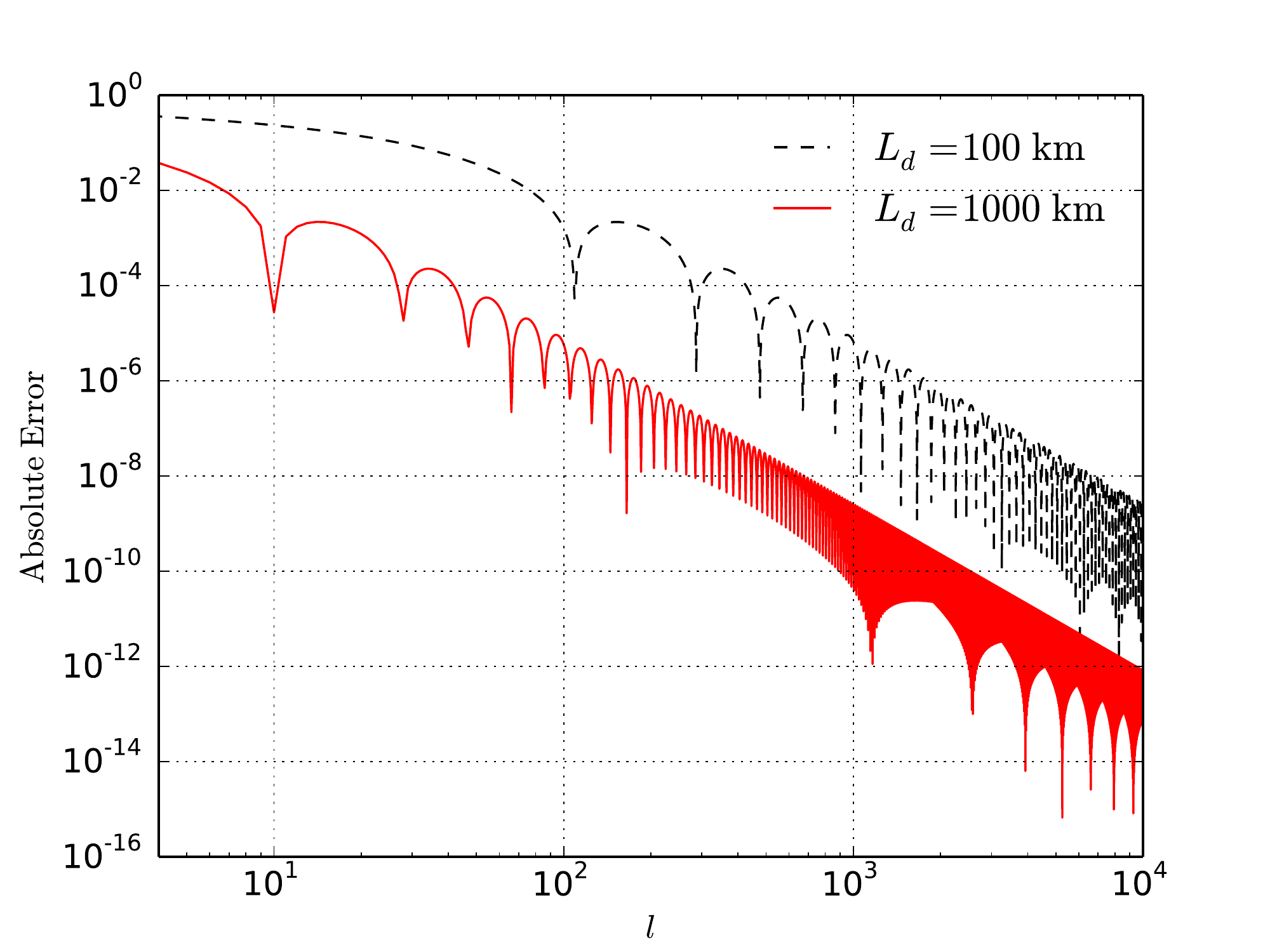}
    \caption{The absolute error (Eq.~\ref{eq:abs-error}) versus the number of terms for $\gamma=\gamma^*$. Curves are generated using the split sum, for $L_d = 1000 \ \text{km}$ and $100$~km  as indicated.}
    \label{figure:AbsoluteError}
\end{figure}

Figure \ref{figure:AbsoluteError} shows the absolute error of the split sum estimate versus the number of terms at $\frac{\gamma}{\gamma^*}=1$, for $L_d=1000 \ \text{km}$ and $L_d=100 \ \text{km}$. The absolute error is computed as 
\begin{equation}
E(l) = \abs{\hat{G}(\gamma,\gamma^*)-G_{l}(\gamma,\gamma^*)},
\label{eq:abs-error}
\end{equation}
where the ``converged'' solution $\hat{G}$ is the value of $G$ obtained after a sufficiently large number of terms, $l'$, is used in the summation, such that the coefficient multiplying the Legendre polynomial in the $l'$th term of the split sum is within quadruple machine precision.  (Actually, it may be shown that for a desired cutoff value of this coefficient, $\epsilon$, the number of terms needed is $l' \simeq \sqrt[3]{ \frac{2}{\epsilon \gamma^{*2}}}$.)  It can be observed that the asymptotic rate of convergence for both values of $L_d$ appears to be similar ($\sim -3.5$ on the log-log plot), though evidently a larger number of terms must be included as $L_d$ decreases.  
Note that $\hat{G}(\gamma,\gamma^*)$ is in close agreement with the value computed using the Numerical Integration Polyalgorithm of Maple$^\text{TM}$ to within 16 decimal points, which is the precision of the Maple integration, as can be seen in Table \ref{Table:splitSumVsMaple}. 
Table \ref{Table:splitSumVsMaple} shows that the split sum approximation matches the numerical integration using the Numerical Integration Polyalgorithm of Maple$^\text{TM}$ over the range $0.001 \leq \frac{\gamma}{\gamma^*} \leq 10$ for $L_d = 1000 \ \text{km}$ ($\gamma^* = 0.15696123$). The agreement improves from 8 significant digits at $\frac{\gamma}{\gamma^*} =0.001$ to 16 significant digits at $\frac{\gamma}{\gamma^*} =10$. For the cases $\cos(\gamma)=0$ and $\cos(\gamma)=-1$, the split sum approximation matches the closed form solutions  (\ref{eqCOS0}) and (\ref{eqCOSM1}), as shown in Table \ref{Table:splitSumVsEQNS} for different values of $L_d$.

\begin{table}[]
\centering
\caption{The integral approximation of the Green's function using the Numerical Integration Polyalgorithm of Maple$^\text{TM}$ and the split sum approximation for $L_d = 1000 \ \text{km}$ ($\gamma^* = 0.15696123$). The stopping criterion used for the split sum is when the absolute value of the factor multiplying the Legendre polynomial reaches quadruple machine precision. Except for the first two entries, the Maple integral was calculated to 16 significant digits. For $\gamma/\gamma^*=0.001$, the maximum number of digits attained was $12$ whereas for $\gamma/\gamma^*=0.005$, it was $15$ digits.}
\label{Table:splitSumVsMaple}
\begin{tabular}{|l|l|l|}
\hline
$\gamma/\gamma^*$ & \textbf{Maple} & \textbf{Split Sum} \\ \hline
0.001 & -1.11851154768 & -1.1185115466790030 \\ \hline
0.005 & -0.862367611582155 & -0.86236761566019138 \\ \hline
0.01 & -0.7520661831497065 & -0.75206618670609104 \\ \hline
0.02 & -0.6418055952187355 & -0.64180559877292642 \\ \hline
0.03 & -0.5773592542889705 & -0.57735925783975817 \\ \hline
0.04 & -0.5316835690653312 & -0.53168357261181198 \\ \hline
0.05 & -0.4963020973589550 & -0.49630209499805700 \\ \hline
0.06 & -0.4674384556430355 & -0.46743845917847821 \\ \hline
0.07 & -0.4430777083432368 & -0.44307770767107485 \\ \hline
0.08 & -0.4220167638214708 & -0.42201676734310684 \\ \hline
0.09 & -0.4034793378064005 & -0.40347933155959809 \\ \hline
0.1 & -0.3869351828355970 & -0.38693518049861692 \\ \hline
0.15 & -0.3237421092103450 & -0.32374210306526552 \\ \hline
0.2 & -0.2796019048871758 & -0.27960190262165774 \\ \hline
0.25 & -0.2459853828209132 & -0.24598538282091331 \\ \hline
0.3 & -0.2190766246998643 & -0.21907661890074209 \\ \hline
0.4 & -0.1780154655314450 & -0.17801546345860769 \\ \hline
0.5 & -0.1477460718583630 & -0.14774607185836305 \\ \hline
0.6 & -0.1243523124660145 & -0.12435230751870802 \\ \hline
0.7 & -0.1057148171283462 & -0.10571481912301701 \\ \hline
0.8 & -9.055025072697948E-002 & -9.0550249089805523E-002 \\ \hline
0.9 & -7.801957852946700E-002 & -7.8019581252887993E-002 \\ \hline
1 & -6.754292534703262E-002 & -6.7542925347032642E-002 \\ \hline
2 & -1.846304821245086E-002 & -1.8463048212450855E-002 \\ \hline
3 & -5.714300085292792E-003 & -5.7143000852927931E-003 \\ \hline
4 & -1.870693766774068E-003 & -1.8706937667740684E-003 \\ \hline
5 & -6.333928838453482E-004 & -6.3339288384534867E-004 \\ \hline
6 & -2.195637338314424E-004 & -2.1956373383144240E-004 \\ \hline
7 & -7.750693244142898E-005 & -7.7506932441429000E-005 \\ \hline
8 & -2.777791606780832E-005 & -2.7777916067808478E-005 \\ \hline
9 & -1.009014887236524E-005 & -1.0090148872365168E-005 \\ \hline
10 & -3.711685976274890E-006 & -3.7116859762750061E-006 \\ \hline
\end{tabular}
\end{table}

\begin{table}[]
\centering
\caption{The split sum approximation for $\gamma=\pi/2$ and $\pi$ compared to Equations (\ref{eqCOS0}) and (\ref{eqCOSM1}).}
\label{Table:splitSumVsEQNS}
\begin{tabular}{|l|l|l|}
\hline
$L_d$ (km) & \textbf{Split Sum ($\cos(\gamma)=0$)} & \textbf{Equation (\ref{eqCOS0})} \\ \hline
300 & -1.4225814713795086E-016 & -1.4225594675545431E-0016 \\ \hline
400 & -6.8995961107067088E-013 & -6.8995960864883886E-0013 \\ \hline
500 & -1.1530544611076218E-010 & -1.1530544610815334E-0010 \\ \hline
600 & -3.5617983195546507E-009 & -3.5617983195574219E-009 \\ \hline
700 & -4.1828668580602192E-008 & -4.1828668580605125E-008 \\ \hline
800 & -2.6799479475630453E-007 & -2.6799479475630768E-007 \\ \hline
900 & -1.1452992927108748E-006 & -1.1452992927108717E-006 \\ \hline
1000 & -3.6839641135260536E-006 & -3.6839641135260568E-006 \\ \hline
1100 & -9.6329029580597044E-006 & -9.6329029580597082E-006 \\ \hline
1200 & -2.1556389233382265E-005 & -2.1556389233382263E-005 \\ \hline
1300 & -4.2781705365952404E-005 & -4.2781705365952408E-005 \\ \hline
1400 & -7.7247341254381796E-005 & -7.7247341254381793E-005 \\ \hline
1500 & -1.2929790801598909E-004 & -1.2929790801598909E-004 \\ \hline
1600 & -2.0346827424618137E-004 & -2.0346827424618136E-004 \\ \hline
1700 & -3.0428682816304554E-004 & -3.0428682816304552E-004 \\ \hline
1800 & -4.3611425145919584E-004 & -4.3611425145919584E-004 \\ \hline
1900 & -6.0302349401310045E-004 & -6.0302349401310041E-004 \\ \hline
2000 & -8.0871962518105109E-004 & -8.0871962518105106E-004 \\ \hline \hline
$L_d$ (km) & \textbf{Split Sum ($\cos(\gamma)=-1$)} & \textbf{Equation (\ref{eqCOSM1})} \\ \hline
600 & -1.6930096452253692E-015 & -1.6890307829549300E-015 \\ \hline 
700 & -1.9949869771970800E-013 & -1.9950267658105915E-013  \\ \hline
800 & -7.1590372634355031E-012 & -7.1590332845768866E-012 \\ \hline 
900 & -1.1609776217297870E-010 & -1.1609776615183131E-010 \\ \hline 
1000 & -1.0797889438705467E-009 & -1.0797889398916860E-009 \\ \hline 
1100 & -6.7027265321956257E-009 & -6.7027265361744720E-009 \\ \hline 
1200 & -3.0722971586707354E-008 & -3.0722971582728503E-008 \\ \hline 
1300 & -1.1152376939660639E-007 & -1.1152376939262755E-007 \\ \hline 
1400 & -3.3703407951899458E-007 & -3.3703407952297341E-007 \\ \hline 
1500 & -8.7963567819438533E-007 & -8.7963567819836426E-007 \\ \hline 
1600 & -2.0379565833152491E-006 & -2.0379565833112706E-006 \\ \hline 
1700 & -4.2803676572173898E-006 & -4.2803676572213691E-006 \\ \hline 
1800 & -8.2844410488393729E-006 & -8.2844410488353935E-006 \\ \hline 
1900 & -1.4967463004072891E-005 & -1.4967463004068912E-005 \\ \hline 
2000 & -2.5504778524850583E-005 & -2.5504778524854560E-005 \\ \hline 
\end{tabular}
\end{table}

\begin{figure}[h]
  \centering
    \includegraphics[width=0.8\textwidth]{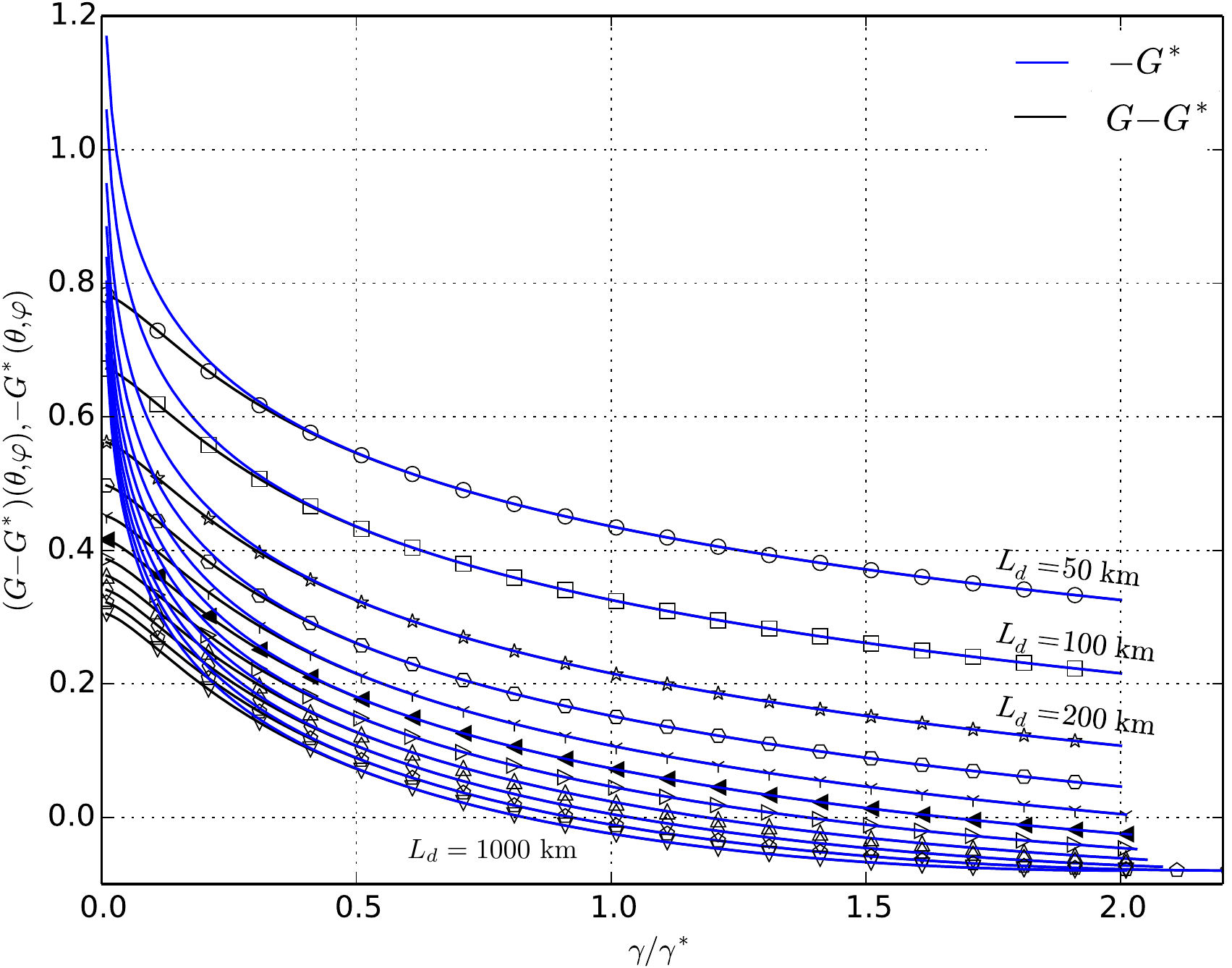}
	\caption{$(G-G^*)$ and $-G^*$ vs. $\gamma/\gamma^*$. Curves are generated using the split sum for $L_d= 50,100,200,300,400,500,600,700,800,900 \ \text{and} \ 1000 \ \text{km}$.}
	\label{fig:GG*}
\end{figure}

Equation (\ref{eq:G-split1}) can be expressed as 
\begin{equation}
G_{l'}(\gamma,\gamma^*) = -\frac{\gamma^{*2}}{4\pi }-\frac{1}{4\pi}\sum_{l=1}^{l'-1} \left[\frac{2l+1}{(l)(l+1)+ \frac{1}{\gamma^{*2}} } - \frac{2l+1}{(l)(l+1)} \right] P_l(\cos(\gamma))  + G^*(\gamma) , 
\label{eq:G-split2}
\end{equation}
where $G^*(\gamma) = \frac{1}{4\pi} \log(\frac{e}{2}(1-\cos(\gamma)))$ is the Green's function of Poisson's equation on the sphere, without the screening term. To explore the departure of $G$ from $G^*$, plots of $(G-G^*)(\gamma, \gamma^*)$ and $-G^*(\gamma, \gamma^*)$ versus $\gamma/\gamma^*$ are presented in Figure \ref{fig:GG*} for the range $50 \text{ km} \le L_d \le 1000 \text{ km}$.

It can be observed from figure \ref{fig:GG*} that for sufficiently large $\gamma/\gamma^*$, the $-G^*$ curve collapses onto the $G-G^*$ curve, indicating that the Green's function of the screened Poisson equation decays at a much faster rate than that of the Poisson equation. This is attributed to the role of the screening term $\psi/L_d^2$ in localizing the solution to a neighborhood of the order $\gamma^*$. This can also be seen in  Figures  \ref{fig:sub1} and \ref{fig:sub2}, where $G$, $G^*$, and $G-G^*$ are plotted versus$\gamma/\gamma^*$ for $L_d=50$ km and $L_d=1000$ km, respectively. It can be observed that $G$ becomes increasingly localized on the sphere as $L_d$ decreases. As such, a compact approximation of $G$ for small values of $L_d$  may prove suitable.

\begin{figure}[hh]
\centering
\begin{subfigure}{0.5\textwidth}
  \centering
  \includegraphics[width=1\linewidth]{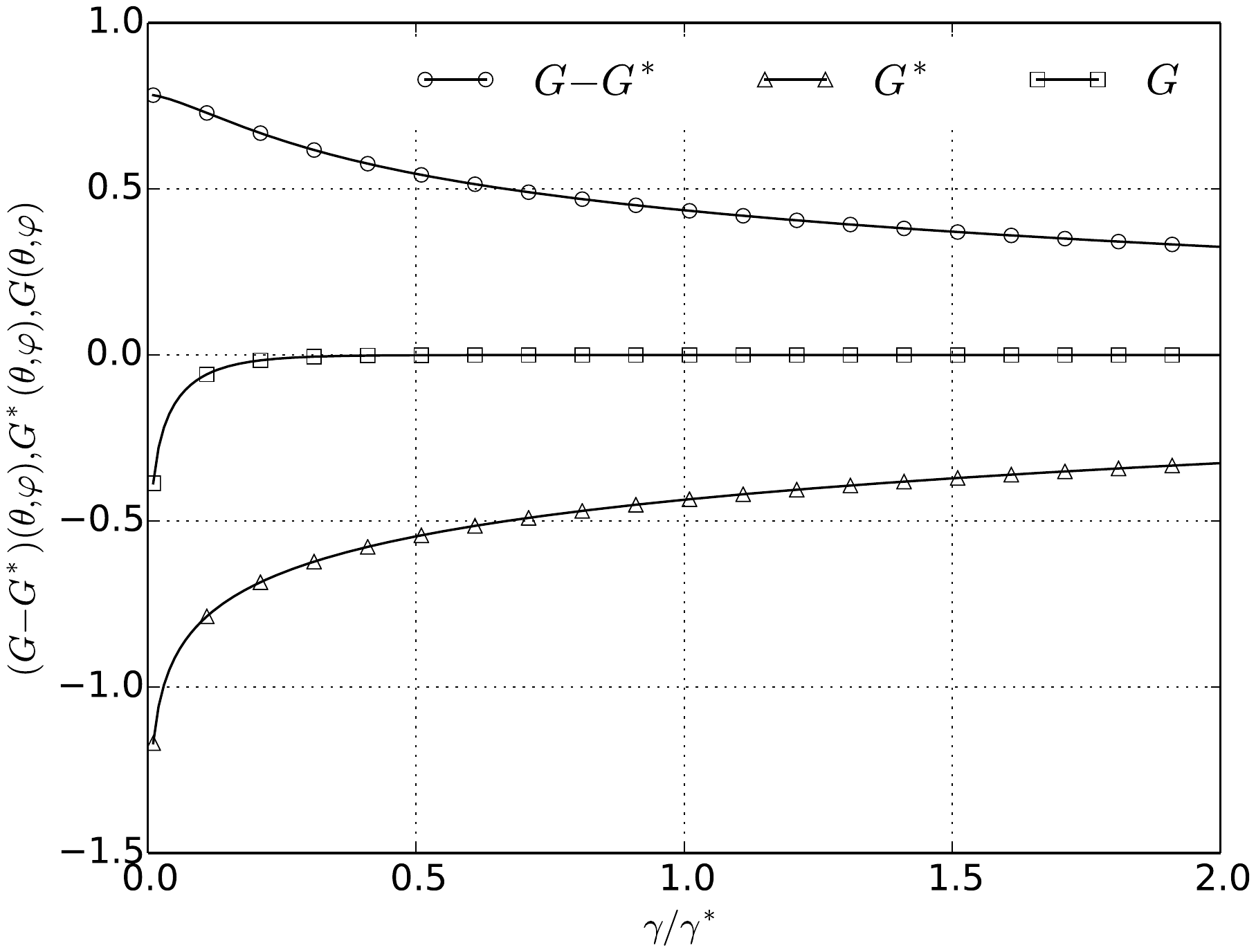}
  \caption{$L_d = 50 \ \text{km}$}
  \label{fig:sub1}
\end{subfigure}%
\begin{subfigure}{0.5\textwidth}
  \centering
  \includegraphics[width=1\linewidth]{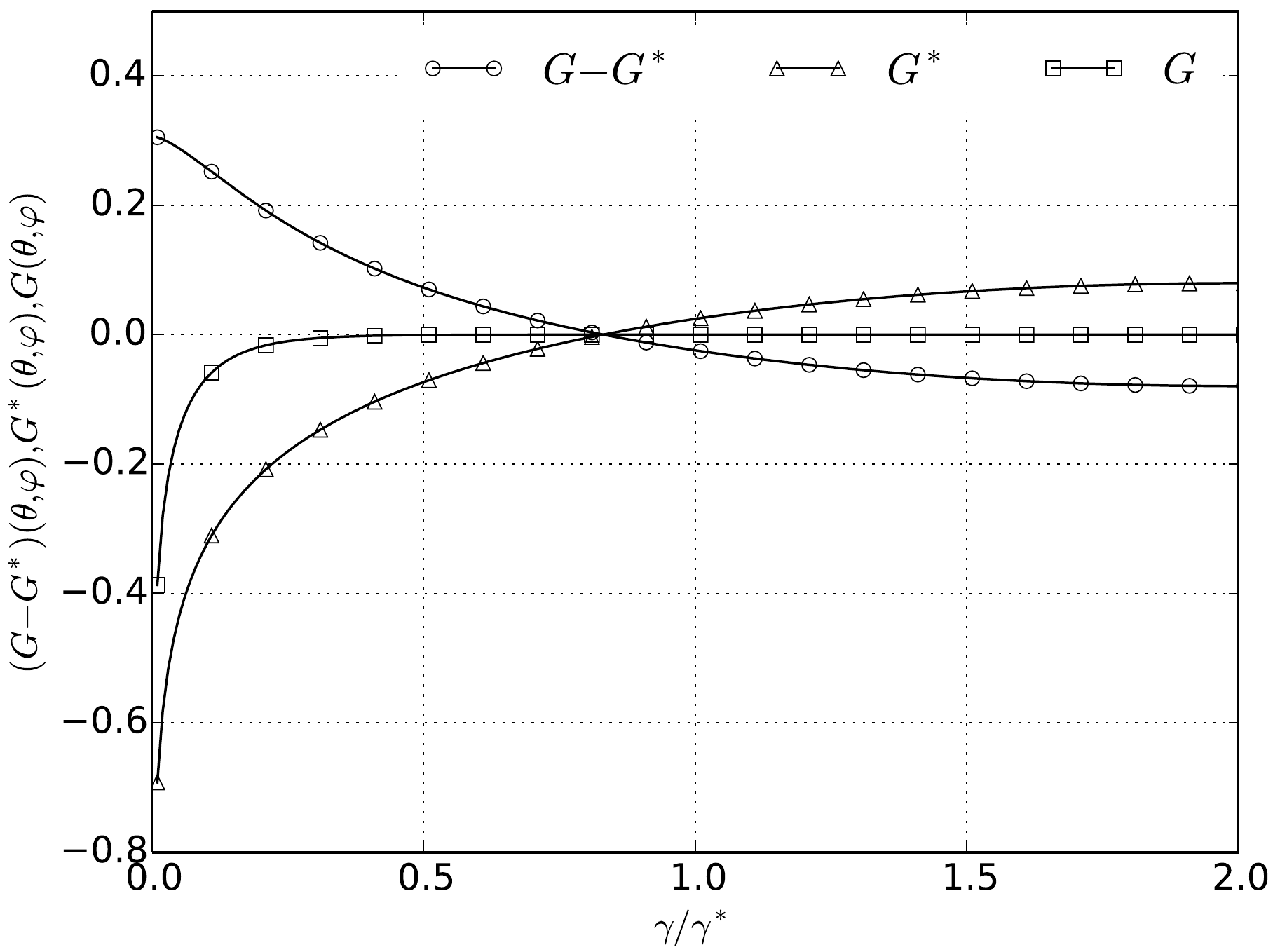}
  \caption{$L_d = 1000 \ \text{km}$}
  \label{fig:sub2}
\end{subfigure}
\caption{$(G-G^*)$, $G^*$ and $G$ vs. $\gamma/\gamma^*$ for (a) $L_d = 50 \ \text{km}$, (b) $L_d = 100 \ \text{km}$. Curves are generated using the split sum.}
\label{fig:Gdecomp}
\end{figure}

%%%%%%%%%%%%%%%%%%%%%%%%%%%%%%%%%%%%%%%%%%%%%%%
%%%%%%%%%%%%%%%%%%%%%%%%%%%%%%%%%%%
%\newpage
\section{Conclusions}
\label{sec:conc}

In this paper, analytical expressions are derived of the Green's function of the screened Poisson's equation on the sphere, 
namely in the form of an integral representation and of a series solution involving Legendre polynomials.  A robust and efficient numerical approximation of the series representation is then developed. This approximation is based
on a splitting of the series representation that is tailored to isolate the singular behavior. Efficiency and robustness of the split 
series approximation was established by showing the rapid decay of the truncation error with the number of terms, and by comparing 
estimates with results obtained using high precision numerical integration. The solutions presented, in both graph and tabular forms, 
for different values of the normalized screening constant versus the normalized angle provide an effective means for accurate evaluation 
of the Green's function.

%\newpage
\section{Acknowledgments}
\label{sec:ackn}

This work is supported by the University Research Board of the American University of Beirut. The authors would like to acknowledge Professor Leila Issa of the Lebanese American University-Beirut for her insightful feedback on the mathematical derivation of the convergence of the Green's Function.

%\newpage
\section{Appendix}
\label{sec:appendix}

In this section, we present tables of the values of the Green's function versus $\gamma/\gamma^*$, for values of the screening constant $L_d$ of $50,100,200,400,800,1000$ km.

% Please add the following required packages to your document preamble:
% \usepackage{longtable}
% Note: It may be necessary to compile the document several times to get a multi-page table to line up properly
\begin{longtable}[c]{|l|l|l|l|}
\caption{$G(\theta,\varphi)$ function of $\gamma/\gamma^*$ for $L_d=50,100,200 \ \text{km}$.The values presented are computed using the split sum. The stopping critera used is when the term contribution of $G(l)$ to the split sum drops down below 1E-20.}
\label{my-label}\\
\hline
$10\gamma/\gamma^*$ & $G(\theta,\varphi)$, $L_d = 50 \ \text{km}$ & $G(\theta,\varphi)$, $L_d = 100 \ \text{km}$ & $G(\theta,\varphi)$, $L_d = 200 \ \text{km}$ \\ \hline
\endfirsthead
\multicolumn{4}{c}%
{{\bfseries Table \thetable\ continued from previous page}} \\
\hline
$10\gamma/\gamma^*$ & $G(\theta,\varphi)$, $L_d = 50 \ \text{km}$ & $G(\theta,\varphi)$, $L_d = 100 \ \text{km}$ & $G(\theta,\varphi)$, $L_d = 200 \ \text{km}$ \\ \hline
\endhead
1 & -0.386281669 & -0.386286557 & -0.386306107 \\ \hline
2 & -0.278953105 & -0.278957963 & -0.278977364 \\ \hline
3 & -0.218435407 & -0.218440190 & -0.218459383 \\ \hline
4 & -0.177384391 & -0.177389115 & -0.177407995 \\ \hline
5 & -0.147127405 & -0.147132024 & -0.147150531 \\ \hline
6 & -0.123747990 & -0.123752505 & -0.123770580 \\ \hline
7 & -0.105126463 & -0.105130859 & -0.105148457 \\ \hline
8 & -8.99792090E-02 & -8.99834707E-02 & -9.00005400E-02 \\ \hline
9 & -7.74669051E-02 & -7.74710327E-02 & -7.74875507E-02 \\ \hline
10 & -6.70094490E-02 & -6.70134276E-02 & -6.70293644E-02 \\ \hline
11 & -5.81887029E-02 & -5.81925362E-02 & -5.82078733E-02 \\ \hline
12 & -5.06933816E-02 & -5.06970622E-02 & -5.07117920E-02 \\ \hline
13 & -4.42856625E-02 & -4.42891903E-02 & -4.43033017E-02 \\ \hline
14 & -3.87800299E-02 & -3.87834013E-02 & -3.87968980E-02 \\ \hline
15 & -3.40292826E-02 & -3.40325013E-02 & -3.40453833E-02 \\ \hline
16 & -2.99149491E-02 & -2.99180150E-02 & -2.99302880E-02 \\ \hline
17 & -2.63405293E-02 & -2.63434462E-02 & -2.63551176E-02 \\ \hline
18 & -2.32266262E-02 & -2.32293960E-02 & -2.32404824E-02 \\ \hline
19 & -2.05073487E-02 & -2.05099769E-02 & -2.05204897E-02 \\ \hline
20 & -1.81276016E-02 & -1.81300901E-02 & -1.81400459E-02 \\ \hline
21 & -1.60410143E-02 & -1.60433669E-02 & -1.60527844E-02 \\ \hline
22 & -1.42083438E-02 & -1.42105669E-02 & -1.42194629E-02 \\ \hline
23 & -1.25962095E-02 & -1.25983069E-02 & -1.26067009E-02 \\ \hline
24 & -1.11760953E-02 & -1.11780716E-02 & -1.11859832E-02 \\ \hline
25 & -9.92354099E-03 & -9.92540177E-03 & -9.93285049E-03 \\ \hline
26 & -8.81749671E-03 & -8.81924666E-03 & -8.82625207E-03 \\ \hline
27 & -7.83978775E-03 & -7.84143247E-03 & -7.84801599E-03 \\ \hline
28 & -6.97468081E-03 & -6.97622448E-03 & -6.98240474E-03 \\ \hline
29 & -6.20851712E-03 & -6.20996533E-03 & -6.21576235E-03 \\ \hline
30 & -5.52941579E-03 & -5.53077273E-03 & -5.53620607E-03 \\ \hline
31 & -4.92701819E-03 & -4.92828898E-03 & -4.93337726E-03 \\ \hline
32 & -4.39227652E-03 & -4.39346535E-03 & -4.39822720E-03 \\ \hline
33 & -3.91727407E-03 & -3.91838606E-03 & -3.92283872E-03 \\ \hline
34 & -3.49507318E-03 & -3.49611230E-03 & -3.50027299E-03 \\ \hline
35 & -3.11958441E-03 & -3.12055484E-03 & -3.12444032E-03 \\ \hline
36 & -2.78545590E-03 & -2.78636138E-03 & -2.78998748E-03 \\ \hline
37 & -2.48797797E-03 & -2.48882244E-03 & -2.49220454E-03 \\ \hline
38 & -2.22300179E-03 & -2.22378876E-03 & -2.22694129E-03 \\ \hline
39 & -1.98686752E-03 & -1.98760070E-03 & -1.99053762E-03 \\ \hline
40 & -1.77634496E-03 & -1.77702762E-03 & -1.77976210E-03 \\ \hline
41 & -1.58857903E-03 & -1.58921431E-03 & -1.59175904E-03 \\ \hline
42 & -1.42104481E-03 & -1.42163562E-03 & -1.42400258E-03 \\ \hline
43 & -1.27150724E-03 & -1.27205648E-03 & -1.27425697E-03 \\ \hline
44 & -1.13798620E-03 & -1.13849656E-03 & -1.14054140E-03 \\ \hline
45 & -1.01872673E-03 & -1.01920078E-03 & -1.02109998E-03 \\ \hline
46 & -9.12171672E-04 & -9.12611722E-04 & -9.14375007E-04 \\ \hline
47 & -8.16938817E-04 & -8.17347143E-04 & -8.18983535E-04 \\ \hline
48 & -7.31800625E-04 & -7.32179382E-04 & -7.33697321E-04 \\ \hline
49 & -6.55665994E-04 & -6.56017161E-04 & -6.57424738E-04 \\ \hline
50 & -5.87564951E-04 & -5.87890449E-04 & -5.89195115E-04 \\ \hline
51 & -5.26634336E-04 & -5.26935910E-04 & -5.28144767E-04 \\ \hline
52 & -4.72106069E-04 & -4.72385378E-04 & -4.73505061E-04 \\ \hline
53 & -4.23296093E-04 & -4.23554680E-04 & -4.24591417E-04 \\ \hline
54 & -3.79595032E-04 & -3.79834353E-04 & -3.80793936E-04 \\ \hline
55 & -3.40459781E-04 & -3.40681203E-04 & -3.41569073E-04 \\ \hline
56 & -3.05406109E-04 & -3.05610913E-04 & -3.06432194E-04 \\ \hline
57 & -2.74002145E-04 & -2.74191494E-04 & -2.74950959E-04 \\ \hline
58 & -2.45862553E-04 & -2.46037584E-04 & -2.46739626E-04 \\ \hline
59 & -2.20643400E-04 & -2.20805145E-04 & -2.21453927E-04 \\ \hline
60 & -1.98037596E-04 & -1.98187015E-04 & -1.98786423E-04 \\ \hline
61 & -1.77770882E-04 & -1.77908878E-04 & -1.78462506E-04 \\ \hline
62 & -1.59598261E-04 & -1.59725663E-04 & -1.60236872E-04 \\ \hline
63 & -1.43300756E-04 & -1.43418350E-04 & -1.43890255E-04 \\ \hline
64 & -1.28682659E-04 & -1.28791173E-04 & -1.29226697E-04 \\ \hline
65 & -1.15568961E-04 & -1.15669085E-04 & -1.16070922E-04 \\ \hline
66 & -1.03803170E-04 & -1.03895516E-04 & -1.04266182E-04 \\ \hline
67 & -9.32452676E-05 & -9.33304182E-05 & -9.36722572E-05 \\ \hline
68 & -8.37699772E-05 & -8.38484848E-05 & -8.41636574E-05 \\ \hline
69 & -7.52651904E-05 & -7.53375498E-05 & -7.56280788E-05 \\ \hline
70 & -6.76305353E-05 & -6.76972195E-05 & -6.79649602E-05 \\ \hline
71 & -6.07761576E-05 & -6.08375885E-05 & -6.10842908E-05 \\ \hline
72 & -5.46215779E-05 & -5.46781630E-05 & -5.49054203E-05 \\ \hline
73 & -4.90947168E-05 & -4.91468236E-05 & -4.93561311E-05 \\ \hline
74 & -4.41309930E-05 & -4.41789707E-05 & -4.43717036E-05 \\ \hline
75 & -3.96725482E-05 & -3.97167169E-05 & -3.98941556E-05 \\ \hline
76 & -3.56675264E-05 & -3.57081772E-05 & -3.58715042E-05 \\ \hline
77 & -3.20694453E-05 & -3.21068510E-05 & -3.22571577E-05 \\ \hline
78 & -2.88366336E-05 & -2.88710471E-05 & -2.90093503E-05 \\ \hline
79 & -2.59317294E-05 & -2.59633835E-05 & -2.60906163E-05 \\ \hline
80 & -2.33212249E-05 & -2.33503379E-05 & -2.34673662E-05 \\ \hline
81 & -2.09750706E-05 & -2.10018407E-05 & -2.11094648E-05 \\ \hline
82 & -1.88663071E-05 & -1.88909180E-05 & -1.89898783E-05 \\ \hline
83 & -1.69707491E-05 & -1.69933737E-05 & -1.70843505E-05 \\ \hline
84 & -1.52666980E-05 & -1.52874909E-05 & -1.53711153E-05 \\ \hline
85 & -1.37346760E-05 & -1.37537836E-05 & -1.38306377E-05 \\ \hline
86 & -1.23572072E-05 & -1.23747632E-05 & -1.24453845E-05 \\ \hline
87 & -1.11186018E-05 & -1.11347299E-05 & -1.11996142E-05 \\ \hline
88 & -1.00047755E-05 & -1.00195894E-05 & -1.00791931E-05 \\ \hline
89 & -9.00308260E-06 & -9.01668682E-06 & -9.07143294E-06 \\ \hline
90 & -8.10216807E-06 & -8.11465998E-06 & -8.16493684E-06 \\ \hline
91 & -7.29183557E-06 & -7.30330430E-06 & -7.34947025E-06 \\ \hline
92 & -6.56292559E-06 & -6.57345345E-06 & -6.61583863E-06 \\ \hline
93 & -5.90721220E-06 & -5.91687558E-06 & -5.95578467E-06 \\ \hline
94 & -5.31730620E-06 & -5.32617469E-06 & -5.36188782E-06 \\ \hline
95 & -4.78656784E-06 & -4.79470600E-06 & -4.82748146E-06 \\ \hline
96 & -4.30903310E-06 & -4.31649960E-06 & -4.34657522E-06 \\ \hline
97 & -3.87934051E-06 & -3.88619037E-06 & -3.91378535E-06 \\ \hline
98 & -3.49267475E-06 & -3.49895777E-06 & -3.52427310E-06 \\ \hline
99 & -3.14470549E-06 & -3.15046805E-06 & -3.17368972E-06 \\ \hline
100 & -2.83154236E-06 & -2.83682698E-06 & -2.85812553E-06 \\ \hline
101 & -2.54968745E-06 & -2.55453324E-06 & -2.57406555E-06 \\ \hline
102 & -2.29599664E-06 & -2.30043929E-06 & -2.31834997E-06 \\ \hline
103 & -2.06764298E-06 & -2.07171570E-06 & -2.08813753E-06 \\ \hline
104 & -1.86208513E-06 & -1.86581826E-06 & -1.88087324E-06 \\ \hline
105 & -1.67703740E-06 & -1.68045904E-06 & -1.69425948E-06 \\ \hline
106 & -1.51044492E-06 & -1.51358063E-06 & -1.52622977E-06 \\ \hline
107 & -1.36045946E-06 & -1.36333290E-06 & -1.37492543E-06 \\ \hline
108 & -1.22541883E-06 & -1.22805159E-06 & -1.23867483E-06 \\ \hline
109 & -1.10382803E-06 & -1.10623989E-06 & -1.11597399E-06 \\ \hline
110 & -9.94342145E-07 & -9.96551535E-07 & -1.00546993E-06 \\ \hline
111 & -8.95751498E-07 & -8.97775294E-07 & -9.05945456E-07 \\ \hline
112 & -8.06967819E-07 & -8.08821312E-07 & -8.16305430E-07 \\ \hline
113 & -7.27011923E-07 & -7.28709381E-07 & -7.35564299E-07 \\ \hline
114 & -6.55002964E-07 & -6.56557290E-07 & -6.62835419E-07 \\ \hline
115 & -5.90148204E-07 & -5.91571393E-07 & -5.97320707E-07 \\ \hline
116 & -5.31734372E-07 & -5.33037337E-07 & -5.38301890E-07 \\ \hline
117 & -4.79119592E-07 & -4.80312451E-07 & -4.85132716E-07 \\ \hline
118 & -4.31726249E-07 & -4.32818126E-07 & -4.37231193E-07 \\ \hline
119 & -3.89034426E-07 & -3.90033819E-07 & -3.94073737E-07 \\ \hline
120 & -3.50576187E-07 & -3.51490826E-07 & -3.55188860E-07 \\ \hline
121 & -3.15930350E-07 & -3.16767370E-07 & -3.20152168E-07 \\ \hline
122 & -2.84717800E-07 & -2.85483736E-07 & -2.88581560E-07 \\ \hline
123 & -2.56597247E-07 & -2.57298069E-07 & -2.60133021E-07 \\ \hline
124 & -2.31261438E-07 & -2.31902618E-07 & -2.34496824E-07 \\ \hline
125 & -2.08433789E-07 & -2.09020371E-07 & -2.11394052E-07 \\ \hline
126 & -1.87865254E-07 & -1.88401842E-07 & -1.90573601E-07 \\ \hline
127 & -1.69331628E-07 & -1.69822428E-07 & -1.71809290E-07 \\ \hline
128 & -1.52630975E-07 & -1.53079881E-07 & -1.54897435E-07 \\ \hline
129 & -1.37581523E-07 & -1.37992075E-07 & -1.39654645E-07 \\ \hline
130 & -1.24019579E-07 & -1.24395001E-07 & -1.25915690E-07 \\ \hline
131 & -1.11797668E-07 & -1.12140981E-07 & -1.13531790E-07 \\ \hline
132 & -1.00783062E-07 & -1.01096965E-07 & -1.02368901E-07 \\ \hline
133 & -9.08561688E-08 & -9.11431712E-08 & -9.23063155E-08 \\ \hline
134 & -8.19092989E-08 & -8.21716881E-08 & -8.32352640E-08 \\ \hline
135 & -7.38454560E-08 & -7.40853068E-08 & -7.50577840E-08 \\ \hline
136 & -6.65772504E-08 & -6.67965026E-08 & -6.76856118E-08 \\ \hline
137 & -6.00259824E-08 & -6.02263910E-08 & -6.10392377E-08 \\ \hline
138 & -5.41207683E-08 & -5.43039427E-08 & -5.50470141E-08 \\ \hline
139 & -4.87977516E-08 & -4.89651484E-08 & -4.96443917E-08 \\ \hline
140 & -4.39993748E-08 & -4.41523511E-08 & -4.47732091E-08 \\ \hline
141 & -3.96738145E-08 & -3.98136066E-08 & -4.03810603E-08 \\ \hline
142 & -3.57743737E-08 & -3.59021080E-08 & -3.64207224E-08 \\ \hline
143 & -3.22589777E-08 & -3.23756879E-08 & -3.28496341E-08 \\ \hline
144 & -2.90897102E-08 & -2.91963431E-08 & -2.96294473E-08 \\ \hline
145 & -2.62324260E-08 & -2.63298432E-08 & -2.67255995E-08 \\ \hline
146 & -2.36563462E-08 & -2.37453328E-08 & -2.41069404E-08 \\ \hline
147 & -2.13337188E-08 & -2.14150138E-08 & -2.17453984E-08 \\ \hline
148 & -1.92395699E-08 & -1.93138252E-08 & -1.96156691E-08 \\ \hline
149 & -1.73513719E-08 & -1.74191932E-08 & -1.76949460E-08 \\ \hline
150 & -1.56488333E-08 & -1.57107714E-08 & -1.59626712E-08 \\ \hline
151 & -1.41136454E-08 & -1.41702143E-08 & -1.44003156E-08 \\ \hline
152 & -1.27293438E-08 & -1.27809985E-08 & -1.29911744E-08 \\ \hline
153 & -1.14810517E-08 & -1.15282264E-08 & -1.17201910E-08 \\ \hline
154 & -1.03553974E-08 & -1.03984670E-08 & -1.05737916E-08 \\ \hline
155 & -9.34028943E-09 & -9.37961975E-09 & -9.53973700E-09 \\ \hline
156 & -8.42486525E-09 & -8.46077342E-09 & -8.60699600E-09 \\ \hline
157 & -7.59931229E-09 & -7.63209584E-09 & -7.76562104E-09 \\ \hline
158 & -6.85479051E-09 & -6.88471902E-09 & -7.00664504E-09 \\ \hline
159 & -6.18333251E-09 & -6.21065332E-09 & -6.32198160E-09 \\ \hline
160 & -5.57775426E-09 & -5.60269475E-09 & -5.70434011E-09 \\ \hline
161 & -5.03158537E-09 & -5.05434494E-09 & -5.14714760E-09 \\ \hline
162 & -4.53897409E-09 & -4.55975213E-09 & -4.64447547E-09 \\ \hline
163 & -4.09467127E-09 & -4.11363610E-09 & -4.19097956E-09 \\ \hline
164 & -3.69393183E-09 & -3.71123487E-09 & -3.78183973E-09 \\ \hline
165 & -3.33246919E-09 & -3.34826056E-09 & -3.41270923E-09 \\ \hline
166 & -3.00642644E-09 & -3.02084047E-09 & -3.07966852E-09 \\ \hline
167 & -2.71233902E-09 & -2.72548850E-09 & -2.77918244E-09 \\ \hline
168 & -2.44705856E-09 & -2.45905718E-09 & -2.50806331E-09 \\ \hline
169 & -2.20776197E-09 & -2.21870988E-09 & -2.26343566E-09 \\ \hline
170 & -1.99189643E-09 & -2.00188843E-09 & -2.04270645E-09 \\ \hline
171 & -1.79717496E-09 & -1.80628723E-09 & -1.84353699E-09 \\ \hline
172 & -1.62151270E-09 & -1.62982572E-09 & -1.66381753E-09 \\ \hline
173 & -1.46304457E-09 & -1.47062817E-09 & -1.50164581E-09 \\ \hline
174 & -1.32008116E-09 & -1.32700273E-09 & -1.35530509E-09 \\ \hline
175 & -1.19111354E-09 & -1.19742449E-09 & -1.22324773E-09 \\ \hline
176 & -1.07475628E-09 & -1.08051634E-09 & -1.10407716E-09 \\ \hline
177 & -9.69787584E-10 & -9.75037384E-10 & -9.96533633E-10 \\ \hline
178 & -8.75078177E-10 & -8.79870343E-10 & -8.99481212E-10 \\ \hline
179 & -7.89633137E-10 & -7.94003419E-10 & -8.11894552E-10 \\ \hline
180 & -7.12542081E-10 & -7.16527948E-10 & -7.32849004E-10 \\ \hline
181 & -6.42990494E-10 & -6.46622145E-10 & -6.61510458E-10 \\ \hline
182 & -5.80230419E-10 & -5.83545989E-10 & -5.97126126E-10 \\ \hline
183 & -5.23609933E-10 & -5.26629684E-10 & -5.39017220E-10 \\ \hline
184 & -4.72515582E-10 & -4.75272321E-10 & -4.86571061E-10 \\ \hline
185 & -4.26419122E-10 & -4.28929697E-10 & -4.39235009E-10 \\ \hline
186 & -3.84819399E-10 & -3.87111510E-10 & -3.96510380E-10 \\ \hline
187 & -3.47288559E-10 & -3.49375445E-10 & -3.57947338E-10 \\ \hline
188 & -3.13416904E-10 & -3.15322934E-10 & -3.23139737E-10 \\ \hline
189 & -2.82855184E-10 & -2.84593071E-10 & -2.91721508E-10 \\ \hline
190 & -2.55280103E-10 & -2.56861588E-10 & -2.63362082E-10 \\ \hline
191 & -2.30390998E-10 & -2.31835107E-10 & -2.37763254E-10 \\ \hline
192 & -2.07933740E-10 & -2.09250714E-10 & -2.14655807E-10 \\ \hline
193 & -1.87667978E-10 & -1.88868629E-10 & -1.93797103E-10 \\ \hline
194 & -1.69382397E-10 & -1.70473691E-10 & -1.74967846E-10 \\ \hline
195 & -1.52878210E-10 & -1.53873220E-10 & -1.57970317E-10 \\ \hline
196 & -1.37981140E-10 & -1.38890732E-10 & -1.42626216E-10 \\ \hline
197 & -1.24539587E-10 & -1.25368757E-10 & -1.28774297E-10 \\ \hline
198 & -1.12408867E-10 & -1.13164270E-10 & -1.16269397E-10 \\ \hline
199 & -1.01460999E-10 & -1.02149539E-10 & -1.04980309E-10 \\ \hline
200 & -9.15805209E-11 & -9.22080745E-11 & -9.47886630E-11 \\ \hline
\end{longtable}
\newpage
% Please add the following required packages to your document preamble:
% \usepackage{longtable}
% Note: It may be necessary to compile the document several times to get a multi-page table to line up properly
\begin{longtable}[c]{|l|l|l|l|}
\caption{$G(\theta,\varphi)$ function of $\gamma/\gamma^*$ for $L_d=400,800,1000 \ \text{km}$.The values presented are computed using the split sum. The stopping critera used is when the term contribution of $G(l)$ to the split sum drops down below 1E-20.}
\label{my-label}\\
\hline
$10\gamma/\gamma^*$ & $G(\theta,\varphi)$, $L_d = 400 \ \text{km}$ & $G(\theta,\varphi)$, $L_d = 800 \ \text{km}$ & $G(\theta,\varphi)$, $L_d = 1000 \ \text{km}$ \\ \hline
\endfirsthead
\multicolumn{4}{c}%
{{\bfseries Table \thetable\ continued from previous page}} \\
\hline
$10\gamma/\gamma^*$ & $G(\theta,\varphi)$, $L_d = 400 \ \text{km}$ & $G(\theta,\varphi)$, $L_d = 800 \ \text{km}$ & $G(\theta,\varphi)$, $L_d = 1000 \ \text{km}$ \\ \hline
\endhead
1 & -0.386384428 & -0.386698574 & -0.386935174 \\ \hline
2 & -0.279055119 & -0.279367000 & -0.279601902 \\ \hline
3 & -0.218536198 & -0.218844444 & -0.219076619 \\ \hline
4 & -0.177483588 & -0.177786931 & -0.178015471 \\ \hline
5 & -0.147224635 & -0.147522002 & -0.147746071 \\ \hline
6 & -0.123842940 & -0.124133401 & -0.124352314 \\ \hline
7 & -0.105218887 & -0.105501644 & -0.105714813 \\ \hline
8 & -9.00688842E-02 & -9.03433040E-02 & -9.05502513E-02 \\ \hline
9 & -7.75536671E-02 & -7.78192431E-02 & -7.80195743E-02 \\ \hline
10 & -6.70931637E-02 & -6.73494935E-02 & -6.75429255E-02 \\ \hline
11 & -5.82692884E-02 & -5.85160889E-02 & -5.87023981E-02 \\ \hline
12 & -5.07707670E-02 & -5.10078520E-02 & -5.11869080E-02 \\ \hline
13 & -4.43598181E-02 & -4.45870832E-02 & -4.47588041E-02 \\ \hline
14 & -3.88509482E-02 & -3.90683748E-02 & -3.92327383E-02 \\ \hline
15 & -3.40969786E-02 & -3.43046039E-02 & -3.44616435E-02 \\ \hline
16 & -2.99794525E-02 & -3.01773753E-02 & -3.03271618E-02 \\ \hline
17 & -2.64018904E-02 & -2.65902579E-02 & -2.67328992E-02 \\ \hline
18 & -2.32849084E-02 & -2.34639104E-02 & -2.35995445E-02 \\ \hline
19 & -2.05626264E-02 & -2.07324829E-02 & -2.08612736E-02 \\ \hline
20 & -1.81799550E-02 & -1.83409173E-02 & -1.84630472E-02 \\ \hline
21 & -1.60905365E-02 & -1.62428729E-02 & -1.63585451E-02 \\ \hline
22 & -1.42551288E-02 & -1.43991308E-02 & -1.45085575E-02 \\ \hline
23 & -1.26403589E-02 & -1.27763264E-02 & -1.28797302E-02 \\ \hline
24 & -1.12177096E-02 & -1.13459527E-02 & -1.14435637E-02 \\ \hline
25 & -9.96272545E-03 & -1.00835599E-02 & -1.01756109E-02 \\ \hline
26 & -8.85435659E-03 & -8.96809902E-03 & -9.05482657E-03 \\ \hline
27 & -7.87442829E-03 & -7.98139721E-03 & -8.06303602E-03 \\ \hline
28 & -7.00720632E-03 & -7.10771605E-03 & -7.18449941E-03 \\ \hline
29 & -6.23903004E-03 & -6.33339258E-03 & -6.40555192E-03 \\ \hline
30 & -5.55801764E-03 & -5.64653799E-03 & -5.71429962E-03 \\ \hline
31 & -4.95380722E-03 & -5.03678387E-03 & -5.10037038E-03 \\ \hline
32 & -4.41734865E-03 & -4.49507311E-03 & -4.55470011E-03 \\ \hline
33 & -3.94072337E-03 & -4.01347689E-03 & -4.06935439E-03 \\ \hline
34 & -3.51698906E-03 & -3.58504499E-03 & -3.63737578E-03 \\ \hline
35 & -3.14005348E-03 & -3.20367469E-03 & -3.25265480E-03 \\ \hline
36 & -2.80456175E-03 & -2.86400085E-03 & -2.90981843E-03 \\ \hline
37 & -2.50580069E-03 & -2.56130029E-03 & -2.60413601E-03 \\ \hline
38 & -2.23961752E-03 & -2.29141000E-03 & -2.33143684E-03 \\ \hline
39 & -2.00234959E-03 & -2.05065659E-03 & -2.08804011E-03 \\ \hline
40 & -1.79076288E-03 & -1.83579559E-03 & -1.87069364E-03 \\ \hline
41 & -1.60199881E-03 & -1.64395850E-03 & -1.67652150E-03 \\ \hline
42 & -1.43352943E-03 & -1.47260714E-03 & -1.50297780E-03 \\ \hline
43 & -1.28311617E-03 & -1.31949317E-03 & -1.34780735E-03 \\ \hline
44 & -1.14877592E-03 & -1.18262402E-03 & -1.20901014E-03 \\ \hline
45 & -1.02875044E-03 & -1.06023205E-03 & -1.08481187E-03 \\ \hline
46 & -9.21479717E-04 & -9.50748275E-04 & -9.73636983E-04 \\ \hline
47 & -8.25578638E-04 & -8.52779020E-04 & -8.74085352E-04 \\ \hline
48 & -7.39816925E-04 & -7.65085628E-04 & -7.84912205E-04 \\ \hline
49 & -6.63100916E-04 & -6.86566520E-04 & -7.05009967E-04 \\ \hline
50 & -5.94457961E-04 & -6.16241363E-04 & -6.33392832E-04 \\ \hline
51 & -5.33022627E-04 & -5.53237507E-04 & -5.69182623E-04 \\ \hline
52 & -4.78024449E-04 & -4.96777589E-04 & -5.11596852E-04 \\ \hline
53 & -4.28777217E-04 & -4.46168735E-04 & -4.59937815E-04 \\ \hline
54 & -3.84669518E-04 & -4.00793273E-04 & -4.13583126E-04 \\ \hline
55 & -3.45156266E-04 & -3.60100210E-04 & -3.71977425E-04 \\ \hline
56 & -3.09751369E-04 & -3.23597807E-04 & -3.34624754E-04 \\ \hline
57 & -2.78021209E-04 & -2.90847180E-04 & -3.01082298E-04 \\ \hline
58 & -2.49578821E-04 & -2.61456298E-04 & -2.70954217E-04 \\ \hline
59 & -2.24078656E-04 & -2.35074913E-04 & -2.43886810E-04 \\ \hline
60 & -2.01212213E-04 & -2.11390012E-04 & -2.19563706E-04 \\ \hline
61 & -1.80703821E-04 & -1.90121791E-04 & -1.97701956E-04 \\ \hline
62 & -1.62307202E-04 & -1.71019958E-04 & -1.78048329E-04 \\ \hline
63 & -1.45802158E-04 & -1.53860659E-04 & -1.60376163E-04 \\ \hline
64 & -1.30991830E-04 & -1.38443531E-04 & -1.44482517E-04 \\ \hline
65 & -1.17700161E-04 & -1.24589249E-04 & -1.30185595E-04 \\ \hline
66 & -1.05769635E-04 & -1.12137255E-04 & -1.17322532E-04 \\ \hline
67 & -9.50593094E-05 & -1.00943726E-04 & -1.05747371E-04 \\ \hline
68 & -8.54430400E-05 & -9.08798320E-05 & -9.53292547E-05 \\ \hline
69 & -7.68078753E-05 & -8.18301341E-05 & -8.59508509E-05 \\ \hline
70 & -6.90527013E-05 & -7.36911607E-05 & -7.75069275E-05 \\ \hline
71 & -6.20869396E-05 & -6.63701576E-05 & -6.99030570E-05 \\ \hline
72 & -5.58294523E-05 & -5.97839353E-05 & -6.30545183E-05 \\ \hline
73 & -5.02075345E-05 & -5.38578897E-05 & -5.68852702E-05 \\ \hline
74 & -4.51560409E-05 & -4.85250930E-05 & -5.13270279E-05 \\ \hline
75 & -4.06165636E-05 & -4.37254967E-05 & -4.63184842E-05 \\ \hline
76 & -3.65367523E-05 & -3.94051967E-05 & -4.18045638E-05 \\ \hline
77 & -3.28696624E-05 & -3.55158154E-05 & -3.77358010E-05 \\ \hline
78 & -2.95731879E-05 & -3.20139188E-05 & -3.40677325E-05 \\ \hline
79 & -2.66095667E-05 & -2.88605006E-05 & -3.07604096E-05 \\ \hline
80 & -2.39449200E-05 & -2.60205252E-05 & -2.77779127E-05 \\ \hline
81 & -2.15488508E-05 & -2.34625259E-05 & -2.50879511E-05 \\ \hline
82 & -1.93940832E-05 & -2.11582283E-05 & -2.26614848E-05 \\ \hline
83 & -1.74561319E-05 & -1.90822248E-05 & -2.04723947E-05 \\ \hline
84 & -1.57130198E-05 & -1.72116779E-05 & -1.84971850E-05 \\ \hline
85 & -1.41450209E-05 & -1.55260641E-05 & -1.67147118E-05 \\ \hline
86 & -1.27344174E-05 & -1.40069278E-05 & -1.51059503E-05 \\ \hline
87 & -1.14653030E-05 & -1.26376763E-05 & -1.36537765E-05 \\ \hline
88 & -1.03233888E-05 & -1.14033865E-05 & -1.23427744E-05 \\ \hline
89 & -9.29584166E-06 & -1.02906370E-05 & -1.11590643E-05 \\ \hline
90 & -8.37113475E-06 & -9.28735335E-06 & -1.00901480E-05 \\ \hline
91 & -7.53890936E-06 & -8.38267351E-06 & -9.12476662E-06 \\ \hline
92 & -6.78985862E-06 & -7.56682175E-06 & -8.25278130E-06 \\ \hline
93 & -6.11561973E-06 & -6.83100325E-06 & -7.46505657E-06 \\ \hline
94 & -5.50867571E-06 & -6.16729949E-06 & -6.75335968E-06 \\ \hline
95 & -4.96227085E-06 & -5.56858458E-06 & -6.11027281E-06 \\ \hline
96 & -4.47033199E-06 & -5.02844159E-06 & -5.52910888E-06 \\ \hline
97 & -4.02739897E-06 & -4.54109431E-06 & -5.00384112E-06 \\ \hline
98 & -3.62856304E-06 & -4.10134044E-06 & -4.52903441E-06 \\ \hline
99 & -3.26941017E-06 & -3.70449447E-06 & -4.09978884E-06 \\ \hline
100 & -2.94597044E-06 & -3.34633683E-06 & -3.71168539E-06 \\ \hline
101 & -2.65467452E-06 & -3.02306603E-06 & -3.36073822E-06 \\ \hline
102 & -2.39231076E-06 & -2.73125738E-06 & -3.04335163E-06 \\ \hline
103 & -2.15599152E-06 & -2.46782497E-06 & -2.75628190E-06 \\ \hline
104 & -1.94311815E-06 & -2.22998847E-06 & -2.49660184E-06 \\ \hline
105 & -1.75135312E-06 & -2.01524199E-06 & -2.26167026E-06 \\ \hline
106 & -1.57859324E-06 & -1.82132658E-06 & -2.04910316E-06 \\ \hline
107 & -1.42294596E-06 & -1.64620656E-06 & -1.85674855E-06 \\ \hline
108 & -1.28270835E-06 & -1.48804668E-06 & -1.68266388E-06 \\ \hline
109 & -1.15634771E-06 & -1.34519212E-06 & -1.52509517E-06 \\ \hline
110 & -1.04248466E-06 & -1.21615074E-06 & -1.38245855E-06 \\ \hline
111 & -9.39877509E-07 & -1.09957716E-06 & -1.25332360E-06 \\ \hline
112 & -8.47408728E-07 & -9.94258130E-07 & -1.13639840E-06 \\ \hline
113 & -7.64072240E-07 & -8.99099007E-07 & -1.03051593E-06 \\ \hline
114 & -6.88962245E-07 & -8.13112649E-07 & -9.34621937E-07 \\ \hline
115 & -6.21263268E-07 & -7.35408548E-07 & -8.47763715E-07 \\ \hline
116 & -5.60241062E-07 & -6.65183222E-07 & -7.69080543E-07 \\ \hline
117 & -5.05234368E-07 & -6.01711747E-07 & -6.97794405E-07 \\ \hline
118 & -4.55647807E-07 & -5.44339969E-07 & -6.33202205E-07 \\ \hline
119 & -4.10945205E-07 & -4.92477454E-07 & -5.74668263E-07 \\ \hline
120 & -3.70643619E-07 & -4.45591553E-07 & -5.21617949E-07 \\ \hline
121 & -3.34308112E-07 & -4.03201284E-07 & -4.73531713E-07 \\ \hline
122 & -3.01546834E-07 & -3.64872591E-07 & -4.29939774E-07 \\ \hline
123 & -2.72006929E-07 & -3.30213538E-07 & -3.90417199E-07 \\ \hline
124 & -2.45370416E-07 & -2.98870390E-07 & -3.54579697E-07 \\ \hline
125 & -2.21350916E-07 & -2.70523657E-07 & -3.22079671E-07 \\ \hline
126 & -1.99690348E-07 & -2.44884887E-07 & -2.92602635E-07 \\ \hline
127 & -1.80156249E-07 & -2.21693597E-07 & -2.65864031E-07 \\ \hline
128 & -1.62539081E-07 & -2.00714510E-07 & -2.41606443E-07 \\ \hline
129 & -1.46650109E-07 & -1.81735132E-07 & -2.19596799E-07 \\ \hline
130 & -1.32319187E-07 & -1.64563517E-07 & -1.99624267E-07 \\ \hline
131 & -1.19393036E-07 & -1.49026292E-07 & -1.81497953E-07 \\ \hline
132 & -1.07733470E-07 & -1.34966811E-07 & -1.65045066E-07 \\ \hline
133 & -9.72159810E-08 & -1.22243534E-07 & -1.50109145E-07 \\ \hline
134 & -8.77283313E-08 & -1.10728593E-07 & -1.36548522E-07 \\ \hline
135 & -7.91693537E-08 & -1.00306416E-07 & -1.24234930E-07 \\ \hline
136 & -7.14478503E-08 & -9.08725823E-08 & -1.13052174E-07 \\ \hline
137 & -6.44816183E-08 & -8.23327042E-08 & -1.02895029E-07 \\ \hline
138 & -5.81965480E-08 & -7.46014734E-08 & -9.36681559E-08 \\ \hline
139 & -5.25258272E-08 & -6.76017820E-08 & -8.52851940E-08 \\ \hline
140 & -4.74092232E-08 & -6.12639113E-08 & -7.76678917E-08 \\ \hline
141 & -4.27924221E-08 & -5.55248469E-08 & -7.07453296E-08 \\ \hline
142 & -3.86264531E-08 & -5.03276070E-08 & -6.44532392E-08 \\ \hline
143 & -3.48671669E-08 & -4.56206770E-08 & -5.87333808E-08 \\ \hline
144 & -3.14747375E-08 & -4.13574632E-08 & -5.35329470E-08 \\ \hline
145 & -2.84132664E-08 & -3.74958411E-08 & -4.88040683E-08 \\ \hline
146 & -2.56503689E-08 & -3.39976971E-08 & -4.45033486E-08 \\ \hline
147 & -2.31568436E-08 & -3.08285770E-08 & -4.05914271E-08 \\ \hline
148 & -2.09063540E-08 & -2.79573058E-08 & -3.70326099E-08 \\ \hline
149 & -1.88751450E-08 & -2.53556820E-08 & -3.37945316E-08 \\ \hline
150 & -1.70417902E-08 & -2.29981953E-08 & -3.08478256E-08 \\ \hline
151 & -1.53869646E-08 & -2.08617656E-08 & -2.81658483E-08 \\ \hline
152 & -1.38932359E-08 & -1.89255083E-08 & -2.57244359E-08 \\ \hline
153 & -1.25448789E-08 & -1.71705263E-08 & -2.35016540E-08 \\ \hline
154 & -1.13277068E-08 & -1.55797206E-08 & -2.14776001E-08 \\ \hline
155 & -1.02289244E-08 & -1.41376129E-08 & -1.96342125E-08 \\ \hline
156 & -9.23698362E-09 & -1.28301973E-08 & -1.79550899E-08 \\ \hline
157 & -8.34146974E-09 & -1.16447962E-08 & -1.64253517E-08 \\ \hline
158 & -7.53298401E-09 & -1.05699325E-08 & -1.50314801E-08 \\ \hline
159 & -6.80304790E-09 & -9.59521707E-09 & -1.37612064E-08 \\ \hline
160 & -6.14400975E-09 & -8.71124506E-09 & -1.26033859E-08 \\ \hline
161 & -5.54896618E-09 & -7.90949883E-09 & -1.15478977E-08 \\ \hline
162 & -5.01168707E-09 & -7.18226900E-09 & -1.05855493E-08 \\ \hline
163 & -4.52655158E-09 & -6.52257182E-09 & -9.70799086E-09 \\ \hline
164 & -4.08848599E-09 & -5.92408433E-09 & -8.90763641E-09 \\ \hline
165 & -3.69291242E-09 & -5.38107958E-09 & -8.17759194E-09 \\ \hline
166 & -3.33569927E-09 & -4.88837015E-09 & -7.51159668E-09 \\ \hline
167 & -3.01311731E-09 & -4.44125758E-09 & -6.90395785E-09 \\ \hline
168 & -2.72180123E-09 & -4.03548572E-09 & -6.34950403E-09 \\ \hline
169 & -2.45871346E-09 & -3.66719832E-09 & -5.84353499E-09 \\ \hline
170 & -2.22111218E-09 & -3.33290218E-09 & -5.38177947E-09 \\ \hline
171 & -2.00652250E-09 & -3.02943159E-09 & -4.96035346E-09 \\ \hline
172 & -1.81271054E-09 & -2.75391820E-09 & -4.57572913E-09 \\ \hline
173 & -1.63766012E-09 & -2.50376275E-09 & -4.22469837E-09 \\ \hline
174 & -1.47955059E-09 & -2.27660957E-09 & -3.90434662E-09 \\ \hline
175 & -1.33673883E-09 & -2.07032347E-09 & -3.61202490E-09 \\ \hline
176 & -1.20774146E-09 & -1.88296911E-09 & -3.34532624E-09 \\ \hline
177 & -1.09121912E-09 & -1.71279202E-09 & -3.10206394E-09 \\ \hline
178 & -9.85962534E-10 & -1.55820190E-09 & -2.88025204E-09 \\ \hline
179 & -8.90880258E-10 & -1.41775669E-09 & -2.67808775E-09 \\ \hline
180 & -8.04986577E-10 & -1.29014888E-09 & -2.49393417E-09 \\ \hline
181 & -7.27391591E-10 & -1.17419319E-09 & -2.32630715E-09 \\ \hline
182 & -6.57291721E-10 & -1.06881426E-09 & -2.17386131E-09 \\ \hline
183 & -5.93961436E-10 & -9.73037095E-10 & -2.03537787E-09 \\ \hline
184 & -5.36745648E-10 & -8.85977292E-10 & -1.90975413E-09 \\ \hline
185 & -4.85052665E-10 & -8.06832823E-10 & -1.79599446E-09 \\ \hline
186 & -4.38348247E-10 & -7.34875882E-10 & -1.69320002E-09 \\ \hline
187 & -3.96150029E-10 & -6.69446276E-10 & -1.60056179E-09 \\ \hline
188 & -3.58022251E-10 & -6.09945150E-10 & -1.51735369E-09 \\ \hline
189 & -3.23571503E-10 & -5.55828938E-10 & -1.44292556E-09 \\ \hline
190 & -2.92442404E-10 & -5.06604481E-10 & -1.37669776E-09 \\ \hline
191 & -2.64314043E-10 & -4.61824246E-10 & -1.31815625E-09 \\ \hline
192 & -2.38896541E-10 & -4.21081975E-10 & -1.26684785E-09 \\ \hline
193 & -2.15928109E-10 & -3.84008964E-10 & -1.22237653E-09 \\ \hline
194 & -1.95172295E-10 & -3.50270563E-10 & -1.18439991E-09 \\ \hline
195 & -1.76415521E-10 & -3.19562904E-10 & -1.15262611E-09 \\ \hline
196 & -1.59464858E-10 & -2.91610125E-10 & -1.12681164E-09 \\ \hline
197 & -1.44146042E-10 & -2.66161732E-10 & -1.10675868E-09 \\ \hline
198 & -1.30301672E-10 & -2.42990295E-10 & -1.09231402E-09 \\ \hline
199 & -1.17789528E-10 & -2.21889243E-10 & -1.08336740E-09 \\ \hline
200 & -1.06481164E-10 & -2.02670977E-10 & -1.07985043E-09 \\ \hline
\end{longtable}

\bibliographystyle{amsplain}
\bibliography{bibfile.bib}

\end{document}